%
%

\magnification=1200

\font\titfont=cmr10 at 14 pt

\font\headfont=cmr10 at 12 pt



\overfullrule=0in

\def\boxit#1{\hbox{\vrule
 \vtop{%
  \vbox{\hrule\kern 2pt %
     \hbox{\kern 2pt #1\kern 2pt}}%
   \kern 2pt \hrule }%
  \vrule}}

  \def\harr#1#2{\ \smash{\mathop{\hbox to .3in{\rightarrowfill}}\limits^{\scriptstyle#1}_{\scriptstyle#2}}\ }

 \def\GG{{{\bf G} \!\!\!\! {\rm l}}\ }

\def\GL{{\rm GL}}

\def\bra#1#2{\langle #1, #2\rangle}

\def\ss{\subset}

\def\half{\hbox{${1\over 2}$}}
\def\smfrac#1#2{\hbox{${#1\over #2}$}}

\def\dim{{\rm dim}}
\def\dist{{\rm dist}}

\def\log{{\rm log}}
\def\Hess{{\rm Hess}}

\def\tr{{\rm tr}}
\def\max{{\rm max}}
\def\min{{\rm min}}

\def\Hom{{\rm Hom\,}}

\def\Sym{{\rm Sym}^2}

\def\Core{{\rm Core}}

\def\rn{\bbr^n}

\def\Int{{\rm Int}}

\def\Symn{{\Sym(\rn)}}

 \def\cd{{\cal C}}

\def\Theorem#1{\medskip\noindent {\bf THEOREM \bf #1.}}
\def\Prop#1{\medskip\noindent {\bf Proposition #1.}}
\def\Cor#1{\medskip\noindent {\bf Corollary #1.}}
\def\Lemma#1{\medskip\noindent {\bf Lemma #1.}}
\def\Remark#1{\medskip\noindent {\bf Remark #1.}}
\def\Note#1{\medskip\noindent {\bf Note #1.}}
\def\Def#1{\medskip\noindent {\bf Definition #1.}}

\def\Ex#1{\medskip\noindent {\bf Example \bf    #1.}}

\def\pf{\medskip\noindent {\bf Proof.}\ }
\def\qed{\hfill  $\vrule width5pt height5pt depth0pt$}
\def\n{\nabla}

   \def\cp{{\cal P}}

   \def\cn{{\cal N}}
\def\cd{{\cal D}}

\def\cp{{\cal P}}
\def\cf{{\cal F}}

\def\vf{\varphi}

\def\wt{\widetilde}
\def\wh{\widehat}

\def\and{\qquad {\rm and} \qquad}

\def\ol{\overline}
\def\bbr{{\bf R}}
\def\bbc{{\bf C}}

\def\a{\alpha}

\def\d{\delta}
\def\e{\epsilon}

\def\g{\gamma}

\def\l{\lambda}

\def\D{\Delta}

\def\G{\Gamma}
\def\O{\Omega}

\def\psh{plurisubharmonic }

\def\lloc{L^1_{\rm loc}}

\def\bo{\partial \Omega}

\def\PSH{{ \rm PSH}}

 \def\ppsh{$\pp$-plurisubharmonic}

\def\Symn{\Sym(\rn)}
 
\def\USC{{\rm USC}}
\def\fa{{\rm\ \  for\ all\ }}

\def\ob{\overline{\O}}

 \def\Span{{\rm Span}}

\def\AA{1}
\def\BB{2}
\def\CC{3}
\def\DD{4}
\def\EE{5}
\def\FF{6}
\def\HH{7}
\def\II{8}
\def\JJ{9}

\def\Cr{C}
\def\CIL{CIL}
\def\HLo{HL$_1$}
\def\HLt{HL$_2$}
\def\HLth{HL$_3$}
\def\HLf{HL$_4$}
\def\HLfi{HL$_5$}
\def\HLs{HL$_6$}
\def\HLse{HL$_7$}
\def\HLe{HL$_8$}
\def\HLn{HL$_9$}

\centerline{\titfont  GEOMETRIC PLURISUBHARMONICITY}
\smallskip
\centerline{\titfont   AND CONVEXITY }
\smallskip

\centerline{\titfont -  AN INTRODUCTION}
\bigskip

\centerline{\titfont F. Reese Harvey and H. Blaine Lawson, Jr.$^*$}
\vglue .9cm
\smallbreak\footnote{}{ $ {} \sp{ *}{\rm Partially}$  supported by
the N.S.F. } 

\vskip .2in

\centerline{\bf ABSTRACT} \medskip
  \font\abstractfont=cmr10 at 10 pt

{{\parindent= .43in
\narrower\abstractfont \noindent

This is an essay on  potential theory for  geometric plurisubharmonic functions. 
It begins with a given closed subset $\GG$ of the Grassmann bundle $G(p,TX)$
of tangent $p$-planes to a riemannian manifold $X$.  This determines a nonlinear
partial differential equation which is convex but never uniformly elliptic ($p<\dim X$).  
A surprising number of results in complex analysis carry over to this more general setting.  
The notions of: a $\GG$-submanifold, an upper semi-continuous 
$\GG$-\psh function, a $\GG$-convex domain,
a $\GG$-harmonic function, and a $\GG$-free submanifold, are defined. 
Results  include a restriction theorem as well as  the existence and uniqueness  of solutions to the Dirichlet Problem for $\GG$-harmonic functions on $\GG$-convex domains.

}}

\vskip.5in

\centerline{\bf TABLE OF CONTENTS} \bigskip

{{\parindent= .1in\narrower\abstractfont \noindent

\qquad \AA.     Introduction    \smallskip

\qquad \BB.    $\GG$-Plurisubharmonicity for Smooth Functions.  \smallskip

\qquad \CC.     $\GG$-Submanifolds and Restriction.  \smallskip

\qquad \DD.   $\GG$-Convexity and the Core.\smallskip

\qquad \EE.   Boundary Convexity.\smallskip

\qquad \FF.   Upper Semi-Continuous $\GG$-Plurisubharmonic Functions.\smallskip

\qquad \HH.   $\GG$-Harmonic Functions and the Dirichlet Problem.\smallskip

\qquad \II.   Geometric Subequations Involving all the Variables.\smallskip

\qquad \JJ.   Distributionally $\GG$-Plurisubharmonic Functions.\smallskip

}}

\vskip .1in

{{\parindent= .3in\narrower

\noindent
{\bf Appendices: }\medskip

A.  Geometric Subequations.
\smallskip

B.  The Linear-Geometric Case.
\smallskip

}}

\vfill\eject


\noindent{\headfont \AA.\  Introduction}
\medskip

  In a recent series of papers [\HLo]--[\HLe] the authors have studied certain   aspects of degenerate
 non-linear elliptic  partial differential equations and ``subequations''. The results 
include the   development of a generalized potential theory, a restriction theorem, 
and solutions to the Dirichlet Problem.  An important special case -- and, in fact,  the  
motivating case -- of all these results is the ``geometric''  one,  in which the equation is determined
by  a distinguished family $\GG$ of tangent $p$-planes on a manifold (as we explain below).
There are many interesting geometric cases coming, for instance, from the theory of 
calibrations, from almost complex and quaternionic geometry, and from $p$-convexity
in riemannian and hermitian geometry. However, these examples will not be 
emphasized here since they occur in profusion in the earlier papers.

One  aim of this paper is to  collect  together
  the various results in the geometric case.
   Because of their importance as motivation
and their usefulness in non-geometric cases,
we thought it would be helpful to present  them in a coordinated fashion.
This exposition also includes several new theorems.

Given an $n$-dimensional riemannian manifold $X$, let $G(p,TX)$ denote the Grassmann bundle
whose fibre at a point $x$ is the set of $p$-dimensional subspaces of the tangent space $T_xX$.
The starting point is to distinguish  a subset $\GG\ss G(p,TX)$ determining the particular ``geometry''.
Then, for example, one defines the {\sl $\GG$-submanifolds} to  simply 
be those $p$-dimensional submanifolds
$M$ of $X$ with $T_xM\in\GG$ for all $x\in M$.  There is also the  analytical notion of a 
{\sl $\GG$-\psh function},  defined for smooth functions $u$ by using the riemannian hessian
$\Hess_x u$.   For each $W\in  G(p,T_xX)$, one can restrict this quadratic form on $T_xX$  
to $W$ and take its trace.  We then  define $u\in \PSH_{\GG}^\infty(X)$, the set of {\sl smooth
$\GG$-\psh functions on $X$}, by requiring that:
$$
\tr_W \Hess_x u\ \geq\ 0 \qquad \forall\, W\in\GG_x, \ \forall\, x\in X.
\eqno{(\AA.1)}
$$
The set $\cp(\GG_x) \ss \Sym(T_xX)$ of $\GG$-positive quadratic forms (i.e., those satisfying (\AA.1))
is a closed convex cone with vertex at the origin but it is never uniformly elliptic, unless  $p=\dim X$.

The smooth theory, i.e., the study of $\PSH_{\GG}^\infty(X)$, is for the most part a straightforward
extension of standard  results in complex analysis -- where $\GG$ is simply the set of complex 
lines in $\bbc^n$,  and the  functions $u\in \PSH_{\GG}^\infty(X)$  are  the standard smooth plurisubharmonics
on a domain $X\ss \bbc^n$.
In Section \DD \ the existence of various kinds of exhaustion functions for $X$
 are characterized in terms of $\GG$-convex hulls and the $\GG$-core. 
 The  $\GG$-core is empty if and only if $X$ admits a smooth strictly 
 $\GG$-\psh function (Definition \DD.1 and Theorem \DD.2).
We recall the   notion of a {\sl $\GG$-free submanifold} which generalizes the
notion of a totally real submanifold in complex analysis.  The maximal possible
dimension of  such submanifolds provides an upper bound on the homotopy type
of strictly $\GG$-convex manifolds (Theorem \DD.16).
  In Section \EE \  the $\GG$-convexity
 of the boundary of a domain is defined and related to the second fundamental form of the 
 boundary, and also to properties of local defining functions for the boundary.

The notion of $\GG$-plurisubharmonicity for a  general upper semi-continuous function $u$ is defined
in Section \FF\    by requiring that each ``viscosity''  test function $\vf$ for $u$ at each point $x\in X$ satisfies
(\AA.1) (cf.  [\Cr],[\CIL]). A key {\sl positivity condition} (Remark \FF.3)
is satisfied, which ensures that smooth $\GG$-\psh functions are also 
$\GG$-\psh in the second sense (cf.  Lemma \FF.2).  A surprising number of the basic properties
of plurisubharmonic functions in complex analysis carry over to the general geometric case, provided that
$\GG$ is a closed set which locally surjects onto $X$ (Theorem \FF.5). 

Under the additional (but still quite weak)
assumption that $\GG$ admits a smooth neighborhood retraction which preserves the 
fibres of the projection $\pi:G(p,TX)\to X$, {\bf restriction holds} in the sense that for any 
upper semi-continuous $u\in\PSH_\GG(X)$ and any minimal $\GG$-submanifold 
$M\ss X$, {\sl the restriction $u\bigr|_M$ is subharmonic for the riemannian Laplacian $\D_M$
on $M$} (Theorem \FF.7).  That is, $u\bigr|_M$ is subharmonic in any of the many (equivalent) classical senses.
For instance, $u\bigr|_M$ is
``sub-the-$\D_M$-harmonics''.
Finally, if each $W\in \GG$ is the tangent space to some minimal $\GG$-submanifold $M$,
then the converse to restriction also holds.  This justifies the terminology ``plurisubharmonic'.

Next we discuss the solution to the Dirichlet problem on domains $\O\ss\ss X$ with smooth strictly
$\GG$-convex boundary and no core.

A smooth function $u$ is $\GG-${\bf harmonic} if in addition to the inequality (\AA.1) holding,  at each
point $x$ there exists a $W\in\GG_x$ such that equality holds, i.e., $\tr_W \Hess_x u=0$.  
In terms of the set 
$\cp(\GG_x)$ defined by  (\AA.1), this is the requirement that 
$\Hess_x u\in \partial \cp(\GG_x)$ at each point $x$.

The notion of the Dirichlet dual $\wt{\cp(\GG)}$ of $\cp(\GG)$, defined in (\HH.1), enables one 
to extend this notion of $\GG$-harmonicity to general continuous functions
since $\partial \cp(\GG) =  \cp(\GG) \cap(-\wt{\cp(\GG)})$ and  $\wt{\cp(\GG)}$ satisfies the positivity condition
required of a subequation (see Section \HH).    
First, we give a proof of the maximum principle for any upper semi-continuous function $u$
which is $\wt{\cp(\GG)}$-subharmonic (much weaker  than $\GG$=\psh) under our hypothesis 
that the $\GG$-core is empty.  This easily established result is a precursor to comparison.
This notion of $\wt{\cp(\GG)}$-subharmonic is referred to as  {\sl  dually $\GG$-\psh} in this paper.

As long as 
$\GG$ is in a weak sense modeled on a euclidean case $\GG_0\ss G(p,\rn)$,
{\bf both existence and uniqueness hold for the Dirichlet Problem for $\GG$-harmonic functions on $\O$} (see 
Definition \HH.5 and Theorem \HH.6). An outline of our proof from [\HLse] is provided in Section \HH.

Since each closed convex set in a vector space $V$ (in our case $\Sym(T_xX)$) is the intersection 
of its supporting closed half-spaces, linear subequations can be made to play a special role in 
understanding our $\GG$-subequations.  This is seen in Sections \II \ and \JJ.  

In Section \II \ we consider the case where each $\GG_x$ {\sl involves all the variables} in the
tangent space $T_xX$. This means there does not exist a proper linear subspace
$W\ss T_xX$ with $\GG_x\ss \Sym(W)$, and it is equivalent (see Lemma \II.1) to the
condition that there exists  $A\in\Span\,\GG$ with $A>0$.  Under the mild condition of 
regularity (Definition \FF.8), this enables one to write the subequation $\cp(\GG)$ locally 
as the intersection of a family of uniformly elliptic subequations (Corollary \II.3), a
fact that has many consequences. One is the Strong Maximum Principle for 
$\GG$-\psh functions (see Theorem \II.5).

There is a distributional notion of $\GG$-plurisubharmonicity  (but not of $\GG$ harmonicity).
In Section \JJ \ we prove that $\GG$-\psh functions and distributionally 
$\GG$-\psh functions are equivalent in a sense made very precise by Theorem \JJ.2
under the hypothesis that   $\GG$  involves all the variables and is regular.
Strict $\GG$-pluri-subharmonicity can also be defined distributionally and is again equivalent 
to the viscosity definition (Theorem \JJ.8). Section \JJ \ concludes with a  local-to-global
result (of Richberg type [R])  for  $C^\infty$ approximation of strictly $\GG$-\psh functions.

Some of the technical issues involving the various hypotheses on $\GG$, 
such as:  $\GG$ closed, $\GG$ locally surjective onto $X$, $\GG$ having a fibre-preserving 
neighborhood retract, or $\GG$ modeled on a euclidean case $\GG_0$, are discussed in Appendix A,
in conjunction with a discussion of the concept of a subequation (Definition A.2)  in the geometric case.

In appendix B we characterize the subequations which are both linear and geometric under
the weak notion of local jet equivalence (Proposition B.4).

   Finally we note that the extreme case, where $\GG = G(p,TX)$  is chosen to be the full grassmann bundle,
is a basic $\GG$-geometry.  There are many additional results specific to this case which are discussed
   in a separate but companion paper [HL$_9$].
    In that paper we use the classical terminology: p-plurisubharmonicity, p-convexity, etc.

  

 \vskip .3in

\noindent{\headfont \BB.\    $\GG$-Plurisubharmonicity for Smooth Functions.}
\medskip
 
 This concept will be  developed  in stages.  We begin with the basic case.
\medskip

\centerline {\bf  Euclidean Space.}
\medskip

Suppose $V$ is an $n$-dimensional real inner product space, and fix an integer $p$, with
$1\leq p\leq n$.  Let $\Sym(V)$ denote the space of symmetric endomorphisms
of $V$.  Using the inner product, this space is identified with the space of quadratic forms
on $V$.   Let $G(p,V)$ denote the set of $p$-dimensional subspaces of $V$.  For $W\in G(p,V)$, the 
 $W$-{\bf trace} of $A$, denoted $\tr_W A$, is the trace  of the restriction $A\bigr|_W$ of $A$
 to $W$.
 
 We identify the Grassmannian $G(p,V)$ with a subset of $\Sym(V)$ by identifying a subspace $W$
 with orthogonal projection $P_W$ onto the subspace $W$.  The natural inner product on $\Sym(V)$
 is defined by using the trace,  namely $\bra AB = \tr(AB)$.  Under this identification we have
 $$
 \tr_W A\ =\ \bra A {P_W}
 \eqno{(\BB.1)}
 $$

  Let $D^2_x u$ denote the second derivative of a function $u$ at $x\in V$.
  
  \Def{\BB.1}  Suppose that $\GG$ is a closed subset of the Grassmannian $G(p,V)$.
  \medskip
  
  (a) A form $A\in\Sym(V)$ is {\bf $\GG$-positive} if 
  $$
 \tr_W A\ \geq \ 0 \quad \forall\, W\in\GG.
 \eqno{(\BB.2)}
 $$
  
    \medskip
  
  (b)  A smooth function $u$ defined on an open subset $X\ss V$ is said to be 
  
  \qquad {\bf  $\GG$-\psh} 
  if
   $$
 \tr_W D^2_x u\ \geq \ 0 \quad \forall\, W\in\GG \ \  {\rm and \ \ } \forall\, x\in X.
 \eqno{(\BB.3)}
 $$
  
  Let $\cp(\GG)$ denote the set of all $\GG$-positive forms $A\in\Sym(V)$, and let $\PSH^\infty_\GG(X)$
  denote the set of all smooth $\GG$-\psh function on $X$. If $\tr_WA>0$ for all $W\in\GG$, then $A$ is said
  to be {\bf $\GG$-strict}.  Similarly, if the inequalities in (\BB.3) are all  strict, then $u$ is said to be 
  {\bf strictly $\GG$-plurisubharmonic}.
  
  \medskip
  
  Note that: 
  \smallskip
  \centerline
  {$u\in \PSH_{\GG}^\infty(X) \quad\iff\quad D_x^2 u \in \cp(\GG) \qquad  \forall\,   x\in X$, \ \ 
  and}
  \smallskip
  \centerline
  {
  $u$ is $\GG$-strict $ \quad\iff\quad D_x^2 u \in \Int \cp(\GG)  \qquad  \forall\,    x\in X$
  \quad }\smallskip

The next result justifies the terminology. We shall say that a function $u$ is 
{\sl subharmonic on an affine subspace ${\bf W}$} if
$\D_{\bf W}\left( u\bigr|_{{\bf W}\cap X} \right) \geq0$ 
where $\D_{\bf W}$ is the euclidean Laplacian on ${\bf W}$.
A $p$-dimensional affine subspace ${\bf W}$ is called an {\sl affine $\GG$-plane}
if its corresponding vector subspace $W$ is a $\GG$-plane.

\Prop{\BB.2} {\sl
A function $u\in C^\infty(X)$ is   $\GG$-plurisubharmonic if and only if the restriction
$u\bigr|_{{\bf W}\cap X}$ is subharmonic for all affine $\GG$-planes ${\bf W}\ss\rn$.
}

\pf This is obvious from Condition  (2) since with $v=u\bigr|_{{\bf W}\cap X}$, we have 
$\tr_{W} D^2 u = \D_{\bf W} v$ on ${\bf W}\cap X$.\qed

\vskip .3in\centerline{\bf Riemannian Manifolds.}\medskip

Suppose $X$ is an $n$-dimensional riemannian manifold. 
Then the euclidean notions above carry over with $V= T_xX$ and the ordinary
second derivative  of a smooth  function  replaced by the {\bf riemannian hessian}.  
Now the set $\GG$ will be an arbitrary closed subset of the Grassmann bundle
$\pi:G(p,TX)\to X$.
For $u\in C^\infty(X)$ this is a well defined section of the bundle $\Sym(TX)$
given on tangent vector fields 
$V,W$ by 
$$
(\Hess \, u)(V,W)\ =\ VWu -(\nabla_V W)u,
\eqno{(\BB.4)}
$$
where $\n$ denotes the Levi-Civita connection.
Note that under composition with a smooth function $\vf:\bbr\to\bbr$, 
$$
\Hess\, \vf(u) \ =\ \vf'(u) \Hess\, u + \vf''(u) \n u\circ\n u
\eqno{(\BB.5)}
$$

\Def{\BB.1$'$}  A smooth function $u$ on $X$ is said to be 
 {\bf $\GG$-plurisubharmonic} if $\Hess_x u$ is $\GG_x$-positive (where $ \GG_x = \GG\cap \pi^{-1}(x)$)
 at each point $x\in X$, i.e.,
 $$
 \tr_W \Hess_x u \geq 0\quad \forall\,  W\in \GG_x  \ \ {\rm and}\ \ \forall\,x\in X.
 \eqno{(\BB.3)'}
 $$

  Again let $\PSH^\infty_\GG(X)$
  denote the set of all smooth $\GG$-\psh functions on $X$, and let $\cp(\GG)$ denote
  the subset of $\Sym(TX)$ with fibres $\cp(\GG_x)$, the set of $\GG_x$-positive elements in 
  $\Sym(T_xX)$.  If the inequalities in (\BB.3)$'$ are all strict at $x$, 
  then we say that $u$ is {\bf strictly $\GG$-\psh at $x$}.

\medskip
\noindent
{\bf Exercise \BB.1. (Convex  Composition Property).}   If $\vf\in C^\infty(\bbr)$ is convex and increasing, then $u\in \PSH^\infty_\GG(X) \ \Rightarrow
\ \vf\circ u \in \PSH^\infty_\GG(X)$.
If, furthermore $\vf$ is strictly   increasing and convex,  then $u$ strictly 
$\GG$-psh \ $\Rightarrow$   $\vf\circ u$  strictly $\GG$-psh.

\medskip
\noindent
{\bf Exercise \BB.2.}  Show that if $u\in C^\infty(X)$ is strictly $\GG$-psh at a point $x\in X$, then
$u$ is strictly $\GG$-psh  in a neighborhood of $x$. (See Claim 1 in the proof of Lemma A.3.)

\medskip
\noindent
{\bf Exercise \BB.3.}  
  Take $X\equiv\bbr$ and let $\GG \ss G(1,TX) = X\times G(1,\bbr)$ be defined  by setting 
  $\GG_x = G(1, \bbr)$ if $x\geq 0$ and $\GG_x=\emptyset$ if $x<0$.  Show that $\cp(\GG) \ss X\times \Sym(\bbr)
  =\bbr^2$ has fibres $\bbr$ if $x<0$ and $\bbr^+ =[0,\infty)$ if $x\geq0$.  In particular, note that $\cp(\GG)$ is not a
  closed set even though $\GG$ is closed.

\vskip .3in 

\noindent{\headfont \CC.\    $\GG$-Submanifolds and Restriction.}
\medskip

The appropriate geometric objects (in a sense dual to the  $\GG$-\psh  functions)
 are the minimal $\GG$-submanifolds.  In the euclidean case this enlarges 
 the family of affine $\GG$-planes used in Proposition \BB.2.

\Def{\CC.1}  If $M$ is a $p$-dimensional submanifold of $X$ with $T_xM\in \GG_x$ for all $x\in M$, then $M$ is said
to be a {\bf $\GG$-submanifold.}
\medskip

Restriction holds as follows.

\Theorem{\CC.2} {\sl
If a  function $u\in C^\infty(X)$ is $\GG$-plurisubharmonic, then the restriction of $u$ to every
 minimal $\GG$-submanifold $M$ is subharmonic in the induced riemannian structure on $M$.
}

\Remark {\CC.3} If $\GG$ is determined by a calibration $\phi$, i.e., $\GG$ consists of the $p$-planes
calibrated by $\phi$ (with the orientation dropped), then $\GG$-submanifolds are automatically minimal.
 Recently, Robles [Ro] has shown that if the calibration is parallel, then this remains true for any
 critical set $\GG$ corresponding to a non-zero critical value of the calibration.
 
\pf
Suppose $M\subset X$ is any $p$-dimensional submanifold, and let $H_M$ denote its mean curvature
vector field.  Then
$$
\D_M\left(u\bigr|_M\right)\ =\ \tr_{{{T}}M} \Hess\, u - H_M u.
$$
In particular, if $M$ is minimal, then
$$
\D_M\left(u\bigr|_M\right)\ =\ \tr_{{{T}}M} \Hess \,u.
\eqno{(\CC.1)}
$$
If $M$ is  a $\GG$-submanifold, then $ \tr_{{{T}}M} \Hess \,u\geq0$ and the result follows.\qed

\Remark {\CC.4} If  for every point $x\in X$ and every $p$-plane $W\in \GG_x$,
 there exists a minimal  submanifold $M$ with $T_xM =W$, then the converse to Theorem \CC.2 is true
 (use the formula (\CC.1)).

 \def\P{\GG}
  \def\ppsh{$\GG$-\psh}
  \def\pcx{$\GG$-convex\ }

\vfill\eject

\noindent{\headfont \DD.\   $\GG$-Convexity and the Core.}
\medskip

We  will answer four  questions concerning the existence of \ppsh functions.
\medskip

(1)\ When does there exist $u\in \PSH^\infty_\P(X)$ which is everywhere strict?

\medskip

(2)\ When does there exist $u\in \PSH^\infty_\P(X)$ which is a proper exhaustion for $X$?

\medskip

(3)\ When does there exist $u\in \PSH^\infty_\P(X)$ which is both strict and an exhaustion?

\medskip

(4)\ When does there exist $u\in \PSH^\infty_\P(X)$ which is an exhaustion and  strict near $\infty$?
\medskip
\noindent
The answers illustrate  some of  the flexibility available in constructing \ppsh functions.

First we characterize those manifolds $X$ which admit a smooth strictly \ppsh function.

\Def{\DD.1.  (The Core)} The {\bf $\P$-core} of $X$ is defined to be the subset
$$
\Core_\P(X)\ =\ \{x\in X  :  {\rm no}\ \  u\in \PSH^\infty_\P (X)\ \ {\rm is\ strict\ at\ \ }x\}.
$$

Note that the core is the intersection over $u\in \PSH^\infty_\P (X)$ of the closed sets where the given
$u$ is not strict, and as such is a closed subset of $X$ (see Exercise \BB.2).

\Theorem{\DD.2}  {\sl 
The manifold $X$ admits a smooth strictly \ppsh function \ \ $\iff$\ \ $\Core_\P(X)=\emptyset$.
In fact, there   exists a function $\psi\in \PSH^{\infty}_{\GG}(X)$ which is $\GG$-strict at each point
$x\notin \Core_\GG(X)$.
}
\pf
The implication $\Rightarrow$ is clear from the definition. For the converse choose an exhaustion
of $X$ by compact subsets $K_1\ss K_2\ss \cdots$. 
Given any sequence of smooth functions $u_j\in C^\infty(X)$ and numbers 
 $\e_j >0, \ j \geq1$ with $\sum\e_j<\infty$, if we choose numbers $\d_j>0$ sufficiently small
 that the semi-norms 
$$
\| w\|_{K,j} \ \equiv\ \sup_{K} \sum_{|\a|\leq j} |D^\a u_j| \ <\ \e_j.
 $$
satisfy
$$
\d_j \|u_j\|_{K_j, j} \ \leq \ \e_j,
$$
then $u=\sum_j\d_j u_j$ converges in the $C^\infty$-topology to $u\in C^\infty(X)$.

If $v$ is $\GG$-strict at a point $x$, then $v$ is $\GG$-strict in a neighborhood of $x$ (Exercise \BB.2).
Therefore, if $K$ is a compact set disjoint from $\Core_\GG(X)$, then we can find $v\in \PSH_{\GG}^{\infty}(X)$
which is $\GG$-strict at each point of $K$.  Hence, we may choose $u_j\in  \PSH_{\GG}^{\infty}(X)$
with $u_j$ strict at each point of  $K_j$ of distance $\geq 1/j$ from $\Core_\GG(X)$. Take $\psi \equiv
\sum \d_j u_j$ as above.\qed

\medskip\noindent
{\bf Remark.}  Essentially the same argument proves that there exists 
$\psi\in \PSH^{\infty}_{\GG}(X)$ such that $\tr_W \Hess\, \psi >0$ for all $\GG$-planes
$W$ which do not lie in the {\sl tangential core} (see [\HLo]).

\Def{\DD.3. (The $\GG$-Convex Hull)}  Given a subset $K\ss X$, the {\bf  \pcx hull of $K$} is the set 
$$
\wh K\ =\ \{x\in X : u(x) \leq \sup_K u \ \  \forall \, u\in \PSH^\infty_\P (X)\}.
$$

Note that $\wh{\wh K} = \wh K$ and that $\wh K$ is closed.

\Theorem{\DD.4. ($\P$-Convexity and Exhaustion)} {\sl
The following three conditions are equivalent.
\medskip

(1)\ If $K\ss\ss X$, then $\wh K \ss\ss X$.
\medskip

(2)\ $X$ admits a smooth \ppsh proper exhaustion function $u$.
\medskip

(3)\ For some neighborhood of $\infty$,  $X-K$ with $K$ compact,

\qquad  there exists 
$u\in  \PSH^\infty_\P(X-K)$ with $\lim_{x\to\infty}u (x)  = + \infty$.
\medskip
}

Condition (3) is a weakening of condition (2) to a local condition at $\infty$ in the one-point compactification
$\overline X = X\cup\{\infty\}$.

\Def{\DD.5} We say that {\bf $X$ is \pcx} if one of the equivalent condition in Theorem \DD.4 holds.

\medskip

The implication (3) \ $\Rightarrow$ (2) is immediate from the next  (stronger) result.
Here $K$ is a compact subset of $X$.

\Lemma{\DD.6}  {\sl Given $v\in   \PSH^\infty_\P(X-K)$ with $\lim_{x\to\infty} v(x) = + \infty$, there exists
$u\in  \PSH^\infty_\P(X)$ such that $u=v$ in a neighborhood of $\infty$.
}

\pf    For  $c$ sufficiently large, $v$ is smooth and  $\P$-\psh outside the compact set    $\{x\in X : v(x)\leq c-1\}$.  Pick a convex increasing function $\varphi\in C^\infty(\bbr)$   with $\varphi \equiv c$
 on a neighborhood of $(-\infty, c-1]$ and $\varphi(t)=t$ on $(c+1, \infty)$.  Then by  Exercise \BB.1, 
the composition   $\varphi\circ v$ is smooth and  $\P$-\psh on all of $X$.  Moreover, $u=v$ outside 
   the compact set  $\{x\in X : v(x)\leq c+1\}$.\qed
   
   \medskip\noindent
   {\bf Proof that (2) $\Rightarrow$ (1).}  If $K$ is compact, 
   then $c=\sup_K u < \infty$,
   and $\wh K$ is contained in the compact set $\{u\leq c\}$.\qed
   \medskip
   
   The implication (1) $\Rightarrow$ (2) is a construction using the next lemma.
   
   \Lemma{\DD.7}  {\sl Suppose $K\ss X$ is compact.  If $x\notin \wh K$, then  there exists 
$u\in  \PSH^\infty_\P(X)$ satisfying:
\medskip 

(a)\ \ $u\equiv 0$  on a neighborhood of $K$,

\medskip 

(b)\ \ $u(x) > 0$,   and

\medskip 

(c)\ \ $u$  is strict at $x$ if $x\notin \Core_\GG(X)$.
}

\pf  Suppose $x\notin \wh K$. Then there exists $v\in  \PSH^\infty_\P(X)$ with $\sup_K v <0< v(x)$.
Pick $\varphi\in C^\infty(\bbr)$ with $\varphi\equiv 0$ on $(-\infty, 0]$ and with $\varphi>0$ and convex increasing on $(0,\infty)$.  Then $u=\varphi\circ v$ satisfies the required conditions.
Furthermore, assume $h\in  \PSH^\infty_\P(X)$ is strict at $x$. Then take  $\ol v = v +\epsilon h$.
For small enough $\epsilon$,  $\sup_K \ol v <0<\ol v(x)$.  If $\varphi$ is also strictly increasing on $(0,\infty)$, then $u=\varphi \circ \ol v$ is strict at $x$.  \qed

\medskip
\noindent
{\bf Proof that  (1) $\Rightarrow$ (2).}  A \ppsh proper exhaustion function on $X$ is constructed as follows.
Choose an exhaustion of $X$ by compact  $\P$-convex subsets 
$K_1\subset K_2\subset K_3\subset \cdots $ with $K_m\subset K^0_{m+1}$ for all $m$.
By Lemma \DD.7  and the compactness of $K_{m+2}-K_{m+1}^0$, there exists a \ppsh function $f_m\geq0$ on $X$ with  $f_m$ identically zero on a neighborhood of $K_m$   and $f_m>0$ on $K_{m+2}-K_{m+1}^0$.
By re-scaling we may assume $f_m>m$ on $K_{m+2}-K_{m+1}^0$. The locally finite sum 
$f=\sum_{m=1}^\infty f_m$ satisfies (2).\qed

\medskip

Next we characterize the existence of a strict exhaustion function.

\Theorem {\DD.8. (Strict $\P$-Convexity)} {\sl
The following conditions are equivalent:
\medskip

(1)\ $\Core_\P(X)=\emptyset$, and if $K\ss\ss X$, then $\wh K \ss\ss X$,
\medskip

(2)\ $X$ admits  a smooth proper exhaustion function which is strictly $\P$-plurisubharmonic.
}

   \medskip\noindent
   {\bf Proof that (1) $\Rightarrow$ (2).}   Since $\Core_\P(X)=\emptyset$, there exists a 
   strictly \ppsh  function $v$ by Proposition \DD.2.  If $u$ is a  \ppsh exhaustion function
   given by Theorem \DD.4,  then $u+e^v$ is a strict exhaustion.\qed

\Def {\DD.9} We say that $X$ is {\bf strictly \pcx} if one of the equivalent conditions of Theorem \DD.8 holds.

\Cor{\DD.10}  {\sl
Suppose that $\Core_\P(X)=\emptyset$.  If $X$ is \pcx, then $X$ is strictly \pcx.
}

\Theorem{\DD.11. (Strict $\P$-Convexity at Infinity)}
{\sl
The following conditions are equivalent:
\medskip

(1)\ $\Core_\P(X)$ is compact, and if $K\ss\ss X$, then $\wh K \ss\ss X$,
\medskip

(2)\ $X$ admits $u\in \PSH_\P(X)$ with $\lim_{x\to\infty} u(x) =\infty$ and $u$ strict outside a 

\ \ \ \ compact  subset.
\medskip

(3)\ $\Core_\P(X)$ is compact, and $X$ admits $u\in \PSH_\P(X-K)$, for some compact 

\ \ \ \ \  subset $K$, with $\lim_{x\to\infty} u(x) =\infty$.
}

   \bigskip\noindent
   {\bf Proof that (3) $\Rightarrow$ (2).}   Apply Lemma \DD.6.

   \medskip\noindent
   {\bf Proof that (2) $\Rightarrow$ (1).} (Straightforward)

   \medskip\noindent
   {\bf Proof that (1) $\Rightarrow$ (3).}  $\Core_\P(X)\equiv K$  is compact $\Rightarrow$ 
    $\Core_\P(X-K)=\emptyset$.

\Def{\DD.12}  We say that $X$ {\bf is strictly $\P$-convex at infinity} if one of the equivalent condition in Theorem \DD.11
holds.
\medskip

Some of the previous results can be summarized as follows.

\Cor{\DD.13}  {\sl
Suppose Core$_\GG(X)=\emptyset$. Then the following are equivalent.
\smallskip

(1)\ \ $X$ is $\GG$-convex.
\smallskip

(2)\ \ $X$ is strictly $\GG$-convex.
\smallskip

(3)\ \ $X$ is strictly $\GG$-convex at infinity.
}
 
 \pf Use Theorems \DD.4 and \DD.11.\qed

\Prop{\DD.14}  {\sl Suppose $(M,\partial M)$ is a compact connected $\GG$-submanifold-with-boundary in $X$.
If $M$ is minimal (stationary), then
\medskip

(1) If $\partial M=\emptyset$, then  $M\ss \Core_\P(X)$.

\medskip

(2) If $\partial M \neq \emptyset$, then  $M\ss \wh {\partial M}$.
}

\pf  Since  the restriction of any $u\in \PSH^\infty_\P(X)$ to $M$ is subharmonic
on $M$, the maximum principle applies to $u\bigr|_M$.\qed

\medskip
 This proposition provides an analogue of the support Lemma 3.2 in [\HLt]:
$$
{\rm If\ } M {\ \rm is\ a\  minimal \  } \GG {\rm \ submanifold,\  then\ } M\ \ss\ \wh{\partial M} \cup \Core_\GG(X).
$$

The existence question for strictly $\P$-convex manifolds has two sides.  We briefly mention these results
from both [\HLo] and  [\HLf].

\Def{\DD.15. ($\GG$-Free)}  A subspace $V\ss T_ X$ is said to be 
{\bf $\P$-free} if there are no $\P$-planes contained in 
$V$.  The maximal dimension of such a free subspace, taken over all points $x\in X$, is called the 
{\bf free dimension} of $\P$ and is denoted freedim$(\P)$. 
A submanifold $M$ of $X$ is {\bf $\P$-free} if $T_xM$ is
$\P$-free for each $x\in M$.
\medskip

Strict $\P$-convexity of $X$ imposes conditions on the topology of $X$.

\Theorem {\DD.16} {\sl
A strictly $\P$-convex manifold has the homotopy type of a CW complex of dimension $\leq$ freedim$(\P)$.
}
\medskip

The free dimension of $G$ is computed in many examples in [\HLo] and summarized
in [\HLf].

On the other hand, the  existence of  many strictly $\P$-convex  manifolds is guaranteed by another 
result (see Theorem 6.6 in [\HLo]).

\Theorem {\DD.17}  {\sl
Suppose $M$ is a $\P$-free submanifold of $X$.  Then $M$ has a fundamental neighborhood 
system in $X$ consisting of strictly $\P$-convex manifolds, each of which has $M$ as a deformation retract.
}

 \vskip .3in

\noindent{\headfont \EE.\   Boundary Convexity}
\medskip

Suppose that $\O\ss X$ is an open connected set with smooth non-empty boundary $\bo$ contained 
in an oriented riemannian manifold. Fix a closed subset $\GG\ss G(p,TX)$.

\Def{\EE.5}  A $p$-plane $W\in \GG_x$ at $x\in \bo$ is called a {\bf tangential $\GG$-plane at $x$}
if $W\ss T_x(\bo)$.

 Denote by $II=II_{\bo}$ the second fundamental form of the boundary with respect to the 
{\bf  inward pointing } normal $n$.  This is a symmetric bilinear form on each tangent space
 $T_x(\bo)$ defined by 
 $$
II(v,w) \ =\ - \bra {\n_v n}{w}  \ =\   \bra { n}{\n_vW}  
 $$
where $W$ is any vector field tangent to $\bo$ with $W_x=w$.

 \Def{\EE.2}   The boundary $\bo$ is  {\bf $\GG$-convex at a point $x$} if $\tr_W  II_x \geq 0$ for all 
 tangential  $\GG$-planes $W$ at $x$.
 If this inequality is strict, then we say that $\bo$ is  {\bf strictly $\GG$-convex at  $x$}.
 \medskip

 \Def{\EE.3. (Local defining functions)} Suppose $\rho$ is a smooth function on a neighborhood
 $B$ of a point $x\in\bo$ with $\bo\cap B =\{\rho=0\}$
 and $\O\cap B=\{\rho<0\}$.  If $d\rho$ is non-zero on $\bo\cap B$, then $\rho $ is called a {\bf local
 defining function for $\bo$}.   
 
 \Lemma {\EE.4}  {\sl
 If $\rho$ is a local defining function for $\bo$, then for all $x\in\bo \cap B$,}
$$
\Hess_x \rho \bigr|_{T_x(\bo)} \ =\ |\n \rho(x)| II_{x}
 $$
 \pf
 Suppose that $e$ is a vector field on $B$ tangent to $\bo$ along $\bo$, and note that
$II(e,e) = \bra n {\n_e e} = - {1\over |\n \rho|} \bra {\n\rho}{\n_e e}$ and
$-\bra {\n\rho}{\n_e e} = - (\n_e e)(\rho) =  e(e\rho) - (\n_e e)(\rho) = (\Hess\,\rho)(e,e)$. 
\qed

\Cor{\EE.5}  {\sl
The boundary  $\bo$ is $\GG$-convex at a point $x$ if and only if 
$$
\tr_W \Hess\, \rho \ \geq\ 0 \quad {\rm for\ all\ } \GG\!-\! {\rm planes\ } W  \ {\rm tangent\ to\ }\bo \ {\rm at\ } x
  \eqno{(\EE.1)}
 $$
 where $\rho$ is a local defining function for $\bo$.  In particular the condition (\EE.1)
  is independent of the choice  of local defining function $\rho$.  Moreover, the boundary
  is strictly $\GG$-convex at a point $x$ if and only if the inequalities in(\EE.1) are all strict, 
  again with independence of the choice of $\rho$.}
 
\Remark{\EE.6}  If $\bo$ is $\GG$-free at a point $x\in\bo$ (see Definition \DD.15), then $\bo$ is automatically
strictly $\GG$-convex at $x$ since there are no tangential $\GG$-planes $W$ to consider.
For example, in the extreme case $p=n$ (the Laplacian subequation) all boundaries $\bo$
are strict at each point since all hyperplanes in $T_x X$ are $\GG$-free.

\Theorem {\EE.7}  {\sl
Suppose that $\bo$ is strictly $\GG$-convex.  Then there exists a global $\GG$-\psh defining function
$\rho\in C^\infty(\ob)$ which is strict on a collar $\{-\e \leq \rho\leq 0\}$. 
If Core$(\O) = \emptyset$, then $\rho$ can be chosen to be strict on all of $\ob$.
}

\Cor{\EE.8}  {\sl
If $\bo$ is strictly $\GG$-convex, then $\O$ is strictly $\GG$-convex at $\infty$; and if 
Core$(\O) = \emptyset$, then $\O$ is strictly $\GG$-convex.
}
\medskip
\noindent
{\bf Proof of Corollary.}  Suppose that $\rho\in C^\infty(\ob)$ is a defining function for $\bo$.
Then $-\log(-\rho)$ is an exhaustion funtion for $\O$.
Since the function $\psi:(-\infty,0) \to (-\infty,\infty)$ defined by $\psi(t )= -\log(-t)$is strictly convex and
increasing, 
$$
\eqalign
{
-\log(-\rho) \   &{\rm is\ strictly\ } \GG\!-\!{\rm \psh \ at\  points\  in\ } \O\   \cr
&{\rm where\ } \rho\ {\rm is\  strictly\  } 
\GG\!-\!{\rm plurisubharmonic.}
}
\eqno{(\EE.2)}
$$
(See Exercise \BB.1.) \qed

\medskip
\noindent
{\bf Proof of Theorem.}  
Start with an arbitrary defining function $\rho\in C^\infty(\ob)$ for $\bo$.
Set $\wt \rho \equiv \rho +{\l \over 2} \rho^2$ with $\l>0$.
Then at points in $\bo$
$$
\Hess\,\wt\rho \ =\ (1+\l\rho)\Hess\,\rho +\l\n \rho\circ \n\rho \ =\ \Hess\,\rho +\l \n \rho\circ \n\rho.
\eqno{(\EE.3)}
$$
We will show that:
$$
\eqalign
{
 {\rm For\ } \l \ {\rm  sufficiently\ }&{\rm large,\ } \wt \rho = \rho +\smfrac \l 2 \rho^2\ \  {\rm is \ strictly \ }  
\GG\!-\!{\rm plurisubharmonic} \ \cr
& {\rm  at\ every\ boundary \  point\ } x\in \bo.
}
\eqno{(\EE.4)}
$$
It then follows that $\wt\rho$ is strictly $\GG$-\psh in a neighborhood of $\bo$ in $X$, and hence on some collar 
$\{-\e\leq\wt \rho\leq 0\}$ with $\e>0$.  Choose $\psi(t)$ convex and increasing with $\psi(t) \equiv -\e$ if $t\leq -\e$,
and $\psi(t)=t$ if $t\geq - {\e\over 2}$. Then $\psi(\wt \rho)$ is $\GG$-\psh on $\O$ and equal to 
$\wt\rho$ on the collar $\{-\e\leq\wt \rho\leq 0\}$, thereby providing the required defining function. 
If Core$(\O)$ is empty, then add the global strictly $\GG$-\psh function, provided by Theorem \DD.2, 
to $\psi(\wt\rho)$.

It remains to prove (\EE.4). Each $p$-plane $V\in G(p,T_x X)$ can be put in canonical form with respect
to $T_x\bo$.  Let $n$ denote a unit normal to $T_x\bo$ in $T_xX$. Choose an orthonormal basis $e_1,...,e_p$
for $V$ such that  $e_2,...,e_p$ is  an orthonormal basis for $V\cap (T_x\bo)$. Then $e=\cos \theta_V n+\sin\theta_V e_1$ defines an angle $\theta_V$ mod $\pi$ and a unit vector $e_1 \in T_x\bo$.  Now by (\EE.3) we have
$$
\tr_V \Hess\,\wt\rho \ =\ \tr_V \Hess\,\rho +\l \cos^2\theta_V |\n \rho|^2.
\eqno{(\EE.5)}
$$

The inequality $|\cos\theta_V|<\d$ defines a fundamental neighborhood system 
for $G(p,T\bo)$ as a subset of the bundle $G(p,TX)\bigr|_{\bo}$. Intersecting with $\GG\bigr|_{\bo}$ 
we see that $\GG\cap G(p, T\bo)$ has a fundamental neighborhood system in $\GG\bigr|_{\bo}$ 
given by $\cn_\d \equiv \{ V\in \GG_x : x\in \bo \ {\rm and\ } |\cos\theta_V|<\d\}$.
Since $\bo$ is strictly $\GG$-convex, there exists $\eta>0$ such that  $\tr_W\Hess\,\rho \geq 2\eta$ for all 
$W\in \GG\cap G(p,T\bo)$.  Hence for $\d$ small, $\tr_V \Hess\,\rho \geq \eta$ for all $V \in\cn_\d$.
Choose a lower bound $-M$ for $\tr_V \Hess\,\rho$ over all $V\in \GG\bigr|_{\bo}$.

Assume $V\in \GG_x$, $x\in \bo$.  For $|\cos\theta_V|<\d$, 
$\tr_V \Hess \wt\rho \geq \eta+ \l \cos^2\theta_V |\n \rho|^2\geq \eta$.
For $|\cos\theta_V|\geq\d$, $\tr_V \Hess\,\wt\rho \geq -M +\l \d^2 |\n\rho|^2$
which is $\geq \eta$ if $\l$ is chosen large.  This proves (\EE.4).\qed

\Remark{\EE.9}  Simple examples show that strict $\GG$-convexity of $\bo$ 
does not imply that every defining function $\rho$ for $\bo$ is strictly $\GG$-\psh at points 
of $\bo$.  However, the exhaustion $-\log(-\rho)$ is always strictly $\GG$-\psh on a small enough
collar of $\bo$.  For the proof of this, compute $\Hess(-\log(-\rho))$ and mimick the proof of 
Theorem \EE.7 on the hypersurfaces $\{\rho=\e\}$ (see the proof of Theorem 5.6 in [\HLo]).

  \Remark{\EE.10. (Signed Distance)}
    Recall that a defining function $\rho$ for $\Omega$ satisfies  $|\nabla \rho| \equiv 1$ in a neighborhood of $\bo$ if and only if $\rho$ is the signed distance to $\bo$ ($<0$ in $ \Omega$ and $>0$ outside of $\Omega$).   In fact any function $\rho$ with $|\nabla \rho| \equiv 1$  in a riemannian manifold is, up to an additive constant,   the distance function to (any) one of its level sets.
  In this case it is easy to see that
      $$
    \Hess \,\rho \ =\ \left( \matrix{0 & 0 \cr 0 &  II}       \right)
   \eqno{(\EE.6)} 
   $$
   where $II$ denotes the second fundamental form of the hypersurface $H=\{\rho =\rho(x)\}$
   with respect to the normal $n=  - \nabla \rho$ and the blocking in (\EE.6) is with respect to the 
   splitting $T_xX =  N_xH \oplus T_x H$.     For example let $\rho (x) =\|x\|\equiv r$ in $\bbr^n$. 
   Then direct calculation shows that $\Hess\, \rho
 = {1\over r}(I-\hat x \circ \hat x)$ where $\hat x = x/r$.
Moreover, 
$$
\Hess(\rho+\l\rho^2) \ =\ \left(\matrix{ 2\l & 0\cr
0& II
}\right)
\eqno{(\EE.7)} 
   $$
simplifying the proof of (\EE.3).
Moreover, setting $\d =-\rho\geq 0$, the actual distance to $\bo$ in $\O$, we have
$$
\Hess(-\log \d) \ =\   {1\over \d}\left(\matrix{ {1\over \d} & 0\cr
0& II
}\right)
\eqno{(\EE.8)} 
   $$
giving an easy proof of Remark \EE.9 for this $\rho$.  Namely, with 
$\d(x) \equiv \dist(x,\bo)$ we have that
$$
\bo \ {\rm strictly\ } \GG\!-\!{\rm convex} \quad\Rightarrow\quad -\log \d
\ {\rm is\ strictly\ } \GG\!-\!{ \rm psh\  in\ a 
\ collar.}
\eqno{(\EE.9)} 
   $$

\Remark{\EE.11. ($\GG$-Parallel)}
If $\GG$ is parallel as a subset of $G(p,TX)\ss \Sym(TX)$, then a weakened form of the converse
to (\EE.9) is true. Namely,
\medskip
\centerline
{\sl
If $-\log\d$ is $\GG$-\psh in a collar, then $\bo$ is $\GG$-convex at each point.
}
\pf
If $\bo$ is not $\GG$-convex at $x\in\bo$, then with $\rho\equiv -\d$, $\tr_W\Hess_x\rho<0$
for some $W\in\GG_x$ tangential to $\bo$ at $x$.  let $\g$ denote the geodesic segment in
$\O$ which emanates orthogonally from $\bo$ at $x$.  Since $\d$ is the distance function to $\bo$,
$\g$ is an integral curve of $\n \d$.  Let $W_y$ denote the parallel translate of $W$ along $\g$  to $y$.
Then $W_y\in\GG_y$ and $(\n \d)_y\perp W_y$.  Therefore by (\EE.8), $\tr_{W_y}\Hess_y (-\log \d) = 
{1\over \d}\tr_{W_y}\Hess_y (\rho) <0$ for $y$ sufficiently close to $x$.
Hence $-\log \d$ is not $\GG$-\psh near $\bo$.\qed
\vskip .3in
\centerline{\bf Local Convexity of a Domain $\O\ss X$}
 \medskip
 
For simplicity assume that $\Core_\GG(X)$ is empty.
Then for each open subset  $Y\ss X$ the three notions of convexity, namely
$\GG$-convexity, strict $\GG$-convexity, and strict $\GG$-convexity at infinity,
are all equivalent.

\Def{\EE.12} 
A domain $\O\ss X$ is {\bf locally $\GG$-convex} if each point
$x\in \partial \O$ has a neighborhood $U$ in $X$ such that $U\cap\O$ is $\GG$-convex.
\medskip

Small balls are $\GG$-convex and the intersection of two $\GG$-convex domains is again
$\GG$-convex.  Therefore:
$$
{\rm If\ } \O {\ \rm is\ } \GG\!-\!{\rm convex,\ then\ } \O\ {\rm is\ locally\ }  \GG\!-\!{\rm convex}.
\eqno{(\EE.10)}
$$

Using terminology from complex analysis, we formulate the ``Levi Problem'':
For which pairs $X,\GG$ does 
$$
 \O {\ \rm locally \ } \GG\!-\!{\rm convex} \qquad \Rightarrow\qquad \O\ {\rm is\ }  \GG\!-\!{\rm convex}?
\eqno{(\EE.10)}
$$
Even in the euclidean case this is  not always true. Here is a counterexample.

\Ex{\EE.13. (Horizontal convexity in $\bbr^2$)}
Take  $\GG=\{\bbr\times \{0\}\}\ss G(1,\bbr^2)$  a singleton
consisting of the $x_1$-axis.  A domain is $\GG$-convex if and only if all of its horizontal slices 
are connected.
Choose $\O\ss\ss \bbr^2$ with the property that $\bo$ contains the interval $[-1,1]$ 
on the $x_1$-axis, the lower half of the circle of radius 3 about the origin, and the points
$(-2,1), (2,1)$.
This can be done with $\O$ locally $\GG$-convex but not globally $\GG$-convex.
In addition, the boundary of $\O$ can be made $\GG$-convex.
\smallskip

By contrast, one of the main results of [\HLn] is the solution to the Levi Problem in euclidean space
in the extreme case $\GG=G(p,\rn)$.

\vfill\eject

\noindent{\headfont \FF.\   Upper Semi-Continuous $\GG$-Plurisubharmonic Functions.}
\medskip

Let $X$ be a riemannian manifold, and assume that $\GG\ss G(p,TX)$ is a closed subset.
Denote by $\USC(X)$ the space of upper semi-continuous
 $[-\infty,\infty)$-valued functions on $X$.
 By a {\bf test function} for $u\in \USC(X)$ at a point $x$ we mean a $C^2$-function $\vf$,
 defined near $x$, such that $u\leq \vf$ near $x$ and $u(x)=\vf(x)$.

\Def{\FF.1}  A function $u\in\USC(X)$ is {\bf $\GG$-plurisubharmonic} if for each
$x\in X$ and each test function $\vf$ for $u$ at $x$, the riemannian hessian
$\Hess_x\vf$  at $x$  satisfies
$$
\tr_W\Hess_x\vf\ \geq\ 0 \quad \forall\, W\in\GG_x
$$
i.e., $\Hess_x \vf\in \cp(\GG_x)$.
 The space of these functions is denoted by
$\PSH_\GG(X)$.
\medskip

This definition is an extension of Definition \BB.1$'$ because of the following.

\Lemma{\FF.2} {\sl
Suppose $u\in C^2(X)$.  Then for a point $x\in X$, the following are equivalent:}
$$
\tr_W \Hess_x \vf \ \geq\ 0 \quad \forall\, W\in \GG_x \ \ {\rm and\  all\ test\ functions\ } \vf\ {\rm for\ } u \ {\rm at}\ x,
\eqno{(\FF.1)}
$$
$$
\tr_W \Hess_x  u \ \geq\ 0 \quad \forall\, W\in \GG_x,
\eqno{(\FF.2)}
$$
\pf Note that (\FF.1) $\Rightarrow$ (\FF.2) because we can take $\vf=u$ in (\FF.1).
Assume (\FF.2) and that $\vf$ is a test function for $u$ at $x$.  
Then $\psi\equiv \vf-u\geq0$ near $x$ and vanishes at $x$.
Hence $x$ is a critical point for $\psi$, and the second derivative 
or hessian of $\psi$ is a well defined non-negative element of $\Sym(T_xX)$, 
independent of any metric. In particular, $\tr_W \Hess_x \psi\geq0$ for all $W\in G(p,T_xX)$.
Since $\Hess_x\vf = \Hess_xu + \Hess_x \psi$, taking the $W$-trace with $W\in\GG_x$, we 
see that (\FF.2) $\Rightarrow$  (\FF.1).\qed

\Remark{\FF.3. (Positivity)} Let $\cp_x\ss\Sym(T_xX)$ denote the subset of non-negative elements.
Replacing $\cp(\GG)\ss\Sym(TX)$ with a  general closed subset $F\ss\Sym(TX)$, the above (standard)
proof shows that (\FF.2) implies (\FF.1), i.e., $\Hess_x u \in F_x \ \ \Rightarrow\ \ \Hess_x \vf \in F_x$, 
provided that $F$ satisfies the {\bf positivity condition}:
$$
F_x + \cp_x \ \ss\ F_x \fa x\in X.
\eqno{(P)}
$$

There are several equivalent ways of stating the condition (\FF.1).  We record one that is 
particularly useful, and refer the reader to Appendix A in [\HLse]  for the proof as well as the statements
of the other conditions.

\Lemma{\FF.4} {\sl
Suppose $u\in\USC(X)$.  Then $u\notin\PSH_\GG(X)$ if and only if $\exists\, x_0\in X, \a>0$, 
and a smooth function $\vf$ defined near $x_0$ satisfying:
$$
\eqalign
{
u-\vf \ &\leq \ -\a|x-x_0|^2 \qquad{\rm near\ \ } x_0\cr
u-\vf \ &= \ 0 \qquad\qquad\qquad  \ \ \ \ {\rm at\ \ } x_0
}
$$
but with $\tr_W \Hess_{x_0} \vf  <0$ for some $W\in \GG_{x_0}$.
}
 \bigskip
 
 \centerline{\bf Elementary Properties}\medskip

Even though $\GG\ss\Sym(TX)$ is closed, the subset 
 $\cp(\GG)\ss\Sym(TX)$ of $\GG$-positive elements may not be closed (see Exercise \BB.3).
 However, by Proposition A.6 below, $\cp(\GG)$ is closed if and only if $\pi\bigr|_\GG$ is a
 local surjection. We make this assumption unless the contrary is stated.
 
 The following basic facts can be found for example in [\HLse, Theorem 2.6].
 In fact they hold with $\cp(\GG)$ replaced by  any subequation (see Definition A.2).

\Theorem{\FF.5}

\smallskip

\item{(a)}  (Maximum  Property)  If $u,v \in \PSH_\GG(X)$, then $w=\max\{u,v\}\in \PSH_\GG(X)$.

\smallskip

\item{(b)}     (Coherence Property) If $u \in \PSH_\GG(X)$ is twice differentiable at $x\in X$, then $\Hess_x u$
is $\GG$-positive.

\smallskip


\item{(c)}  (Decreasing Sequence Property)  If $\{ u_j \}$ is a 
decreasing ($u_j\geq u_{j+1}$) sequence of \ \ functions with all $u_j \in \PSH_\GG(X)$,
then the limit $u=\lim_{j\to\infty}u_j \in \PSH_\GG(X)$.

\smallskip

\item{(d)}  (Uniform Limit Property) Suppose  $\{ u_j \} \ss \PSH_\GG(X)$ is a 
sequence which converges to $u$  uniformly on compact subsets to $X$, then $u \in \PSH_\GG(X)$.

\smallskip

\item{(e)}  (Families Locally Bounded Above)  Suppose $\cf\subset \PSH_\GG(X)$ is a family of 
functions which are locally uniformly bounded above.  Then the upper semicontinuous
regularization $v^*$ of the upper envelope 
$$
v(x)\ =\ \sup_{u\in \cf} u(x)
$$
belongs to $\PSH_\GG(X)$.
 \medskip
  
 \Ex{\FF.6} The following examples show that Properties (c), (d) and (e) require that the set
 $\cp(\GG)$ be closed. Let $X=\bbr$ and $\GG_x =\{T_x\bbr\}\in G(1,TX)$ if $x\geq0$ and $\GG_x=\emptyset$
 for $x<0$.  Note that $\GG$ is a closed set.  Then $\cp(\GG_x) = \Sym(T_xX)\cong \bbr$ for $x<0$ and 
 $\cp(\GG_x) = \{A\in  \Sym(T_xX) : A\geq0\}$ for $x\geq0$. Note that $\cp(\GG)$ is not closed in $\bbr\times \bbr$.
 This subequation is simply  the requirement that $$u''(x)\geq 0 \fa x\geq0.$$
  
 Fix a constant $a>0$ and set 
 $$
 u(x) \ =\   \cases{ 
\ \  \quad 0  \ \qquad {\rm if}\ \ x\geq0, \cr
 x(a-x) \quad {\rm if} \ \ x\leq 0.
 }
 $$
This function fails to be $\GG$-\psh at 0. To see this note that $\vf(x) = x(a-x)$
is a test function for $u$ at $0$ and $\vf''(0)<0$.

 For each $\d>0$ set $v_\d(x) = u(x+\d)+\d$. Note that graph$(v_\d) = {\rm graph}(u)+(-\d,\d)$.
 Then each $v_\d$ is $\GG$-\psh and $v_\d\downarrow u$ as $\d\to 0$. Hence condition (c)  fails.

Now for each $\e>0$, define $u_\e \equiv \min\{u, -\e\}$. Then $u_\e$ is $\GG$-\psh for all $\e$
 and $u_\e\uparrow u$ as $\e \to 0$.  Hence conditions (d) and (e) also fail.

 \bigskip

\centerline{\bf Restriction}
\medskip

Throughout this subsection we assume that $\GG\ss G(p,TX)$  is a closed set 
admitting a smooth neighborhood retraction preserving
 the fibres of the projection $\pi: G(p,TX)\to X$.
 The terminology $\GG$-\psh for $u\in \USC(X)$ is justified by the next result, which extends Theorem \CC.2.

 \Theorem{\FF.7} {\sl
 If $u\in\PSH_\GG(X)$, then for every minimal $\GG$-submanifold $M$, the restriction $u\bigr|_M$
 is $\D$-subharmonic where $\D$ is the Laplace-Beltrami operator in the induced riemannian metric on $M$.
 }
 
 \medskip
 
 This result can be extended to submanifolds $M$ of dimension larger that $p$.
 Let $\GG_M \equiv \{W\in G(p,TM) :  W\in\GG\}$ denote the set of {\bf tangential $\GG$-planes to $M$}.
 This set $\GG_M$ defines a notion of $\GG_M$-plurisubharmonicity for functions $w\in \USC(M)$.

 \Def{\FF.8}  We say that $\GG$ is {\bf regular} if at every point $x_0\in X$,   each element $W_0\in \GG_{x_0}$
 has a local  smooth  extension  to a  section $W(x)$  of $\GG$.

 \Def{\FF.9}   A submanifold $M$ of $X$ is {\bf $\GG$-flat} if the second fundamental form $B$
 of $M$ satisfies
 $$
 \tr\left(B\bigr|_W\right) \ =\ 0\fa {\rm tangential \ } \GG\ {\rm planes}\ W\in \GG_M
 \eqno{(\FF.3)}
 $$

 \Theorem{\FF.10} {\sl
 Suppose  $M$ is a $\GG$-flat submanifold of $X$ and that the subset $\GG_M\ss G(p,TM)$ 
 is regular on $M$.  If $u\in \PSH_\GG(X)$, then $u\bigr|_M \in \PSH_{\GG_M}(M)$.
 }
 \medskip
 
 See Section 8 of [\HLe] for a more complete discussion, including Example 8.4, which shows
 that $\GG_M$ being regular is necessary in Theorem \FF.10.  The proof uses Lemma 8.3 in [\HLe]
 which is stated in this paper as Proposition \II.4 below.

 \vskip .3in

\noindent{\headfont \HH.\     $\GG$-Harmonic Functions and the Dirichlet Problem.}
\medskip

In this section we discuss the Dirichlet problem for {\sl extremal} or {\sl  $\GG$-harmonic} functions.
These are natural generalizations of solutions of the classical homogeneous Monge-Amp\`ere problem,
in both the real and complex cases (and constitute a very special case of the general $F$-harmonic
functions treated in [\HLse]).
To do this we must introduce the {\sl Dirichlet dual}.

\bigskip
\centerline{\bf Dually $\GG$-Plurisubharmonic Functions}
\medskip

We first define the  {\sl Dirichlet dual}  of the subset  $F\equiv \cp(\GG) \ss\Sym(TX)$, 
to be the subset  $\wt F \equiv \wt {\cp(\GG)} \ss\Sym(TX)$ whose fibres are given by
$$
\wt{F_x} \ =\ -(\sim \Int F_x)\ =\ \sim (-\Int F_x).
\eqno{(\HH.1)}
$$
Since 
$$
A\in \Int F_x\qquad\iff\qquad\tr_W A >0 \fa W\in \GG_x,
\eqno{(\HH.2)}
$$
it is easy to see that 
$$
\ \ \ A\in\wt{F_x}\qquad\iff\qquad\tr_W A \geq 0 \quad {\rm for\ some\ \ } W\in \GG_x,
\eqno{(\HH.3)}
$$

 \Def{\HH.1} A smooth function $u$ on $X$ is said to be {\bf dually $\GG$-\psh} if at each point $x\in X$ 
$$
\exists \, W\in \GG_x \  \ {\rm with \ }\   \tr_W \Hess_x u \ \geq\ 0, \qquad {\rm or\ equivalently \ }\ \ \ \ \  \Hess_x u\in\wt{\cp(\GG)}.
$$
More generally a function $u\in \USC(X)$ is {\bf  dually $\GG$-\psh} if for each point $x\in X$ 
and each test function $\vf$ for $u$ at $x$,  
$$
\exists \, W\in \GG_x \  \ {\rm with \ }\   \tr_W \Hess_x \vf \ \geq\ 0, \qquad {\rm or\ equivalently \ }\ \ \ \ \  \Hess_x \vf \in\wt{\cp(\GG)}.
$$
The set of all such  functions is denoted $\wt{\PSH}_\GG(X)$.
\medskip

First note that $\wt{\cp(\GG)}$ satisfies the positivity condition (P), so that as noted in 
Remark \FF.3, if a smooth function $u$ satisfies $\Hess_x u \in \wt{\cp(\GG)}$,
then for each test function $\vf$ for $u$ at $x$, we have $\Hess_x \vf \in \wt{\cp(\GG)}$,
making the second definition an extension of the first definition.  Second,
assuming that $\pi\bigr|_\GG$ is a local surjection as in Definition A.5, it then follows that
not only $\cp(\GG)$, but also  $\wt{\cp(\GG)}$ is closed.  As a consequence, 
$$
{\rm the\  set\ } \ \wt{\PSH}_\GG(X) \  \ {\rm  satisfies 
\ all \ of \ the \  properties \  given \ in\  Theorem  \  \FF.5.}
\eqno{(\HH.4)}
$$
In fact $\wt{\cp(\GG)}$ is a subequation (Definition A.2).

By Theorem \DD.2 if $\Core_\GG(X)=\emptyset$, then $X$ admits a smooth function 
$\psi$ which is strictly $\GG$-\psh at each point.  Of course, $\cp(\GG) \ss\wt{\cp(\GG)}$,
so that the dually $\GG$-\psh functions on $X$ constitute a much larger class than
the $\GG$-\psh functions.  Again we assume that $\pi\bigr|_\GG$ is a local surjection.

\Theorem{\HH.2. (The Maximum Principle for Dually $\GG$-Plurisubharmonic Functions)} {\sl
Suppose $\Core_\GG(X)=\emptyset$.  Then for each compact subset $K\ss X$ and each 
$u\in \wt{\PSH}_\GG(K) \equiv \USC(X) \cap \wt{\PSH}_\GG(\Int K)$ we have:
$$
\sup_K u \ \leq\ \sup_{\partial K} u.
$$}

\medskip
The proof is classical and completely elementary.  Moreover, one can easily see that this maximum principle
is equivalent to the special case of comparison (Theorem \HH.7 below) where $u$ is smooth.

\pf Suppose it fails.  Then there exist a compact set $K$, a function $u\in \wt{\PSH}_\GG(K)$
and a point  $\bar x\in \Int K$ with $u(\bar x) >\sup_{\partial K} u$.  Let $\psi$ be a
smooth strictly $\GG$-psh function  on $X$.  Then for $\e>0$ sufficiently small, the function
$u+\e \psi$ will also have a maximum at some point $x\in\Int K$. Thus $-\e\psi$ is a test function for
$u$ at $x$, and therefore 
$\Hess_x (-\e \psi) \in \wt{\cp_x(\GG)} = -(\sim \Int \cp_x(\GG))$,
i.e., $\Hess_x(\psi) \notin \Int \cp_x(\GG)$
contradicting the strictness of $\psi$ at $x$.\qed

The Convex-Increasing Composition Property
in Exercise \BB.1 not only extends to the upper semi-continuous case,
but also to the much larger class of dually $\GG$-\psh functions.

\Lemma{\HH.3. (Composition Property)} {\sl
Suppose $\vf : \bbr\to\bbr$ is both convex and increasing (i.e., non-decreasing).  Then
$$
u\in \wt{\PSH_{\GG}}(X)
\qquad\Rightarrow\qquad
\vf\circ u\in \wt{\PSH_{\GG}}(X)
\eqno{(a)}
$$
If $\vf$ is also strictly increasing, then in addition to (a) we have that
$$
u\ \ {\rm is\ } \GG\ {\rm strict}
\qquad\Rightarrow\qquad
\vf\circ u\ \  {\rm is\ } \GG\ {\rm strict}
\eqno{(b)}
$$
where we refer ahead to Definition \HH.7 for the notion of strictness.
}
\pf
We can assume that $\vf$ is smooth since it can be approximated by a decreasing
sequence $\vf_\e$ via convolution.  Observe now that:

\centerline {\sl
$\psi$ is a test function for $u$ at $x$ \quad$\iff$\quad $\vf\circ \psi$ is a test function for $\vf \circ u$
at $x$.} 

\noindent
This reduces the proof to the case where $\vf$ and $u$ are both smooth, and formula (\BB.5) applies
with both coefficients $\vf'(u(x))$ and $\vf''(u(x))$ $\geq0$.\qed

\bigskip
\centerline{\bf $\GG$-Harmonics}
\medskip

To understand the next definition note that 
$$
\partial \cp(\GG) \ =\ \cp(\GG)\cap (-\wt{ \cp (\GG)})
\eqno{(\HH.5)}
$$

\Def{\HH.4}  A function $u$ on $X$ is said to be {\bf $\GG$-harmonic} if 
$$
u\in {\PSH}_\GG(X) \and -u \in \wt{\PSH}_\GG(X).
$$

\medskip By (\HH.5)   we see that a $C^2$-function $u$ on $X$ is $\GG$-harmonic if and only if 
$$
\Hess_x u\in \partial \cp(\GG_x)  \fa x\in X.
$$
  
In order to solve the Dirichlet Problem for $\GG$-harmonic functions on domains $\O\ss X$,
we restrict $\GG\ss G(p,TX)$ to be  modeled on a ``constant coefficient''
 case $\GG_0 \ss G(p,\rn)$.

\Def{\HH.5}  A closed subset $\GG\ss G(p,TX)$ is   {\bf   locally trivial 
with fibre} $\GG_0\ss G(p, \rn)$, if in a neighborhood
each point $x\in X$ there exists a local tangent frame field so that under the associated 
trivialization $\phi:G(p, TU)\harr{\cong} {} U\times G(p,\rn)$ we have 
$$
\phi: \GG\bigr|_U \ \harr{\cong} {} \ U\times \GG_0.
$$
This can be formulated somewhat differently.  Let ${\rm Aut}(\GG_0) = \{g\in \GL_n: g(\GG_0)=\GG_0\}$.  Then 
given a  closed subset $\GG\ss G(p,TX)$ which is  locally trivial 
with fibre $\GG_0$, the local tangent frame fields in Definition \HH.5 provide $X$ with a 
 topological ${\rm Aut}(\GG_0)$-structure (see \S  5 in [\HLse]).
 Conversely, if $X$ admits a topological ${\rm Aut}(\GG_0)$-structure, then the euclidean model
 $\GG_0\ss G(p, \rn)$ determines a canonical  closed subset $\GG\ss G(p,TX)$ which is locally trivial
 with fibre $\GG_0$. In other words, a euclidean model can be transplanted to any manifold
 with a topological  ${\rm Aut}(\GG_0)$-structure (again see \S  5 in [\HLse]).

 In the language of [\HLse, \S 6]: ``$\GG$ is locally trivial with fibre $\GG_0$'' means that the subequation $\cp(\GG)$ is locally {\sl jet equivalent} to the constant coefficient subequation $\cp(\GG_0)$.
 
 \medskip

In the next two theorems
$X$ is a riemannian manifold and  $\GG\ss G(p,TX)$  is a closed, locally trivial set with non-empty fibre. 

\Theorem{\HH.6. (The Dirichlet Problem)}  {\sl
Suppose  that $\O\ss\ss X$ is a domain with a
  smooth, strictly $\GG$-convex  boundary $\bo$ and Core$_\GG(\O)=\emptyset$.
 Then the Dirichlet problem for $\GG$-harmonic functions is uniquely solvable on $\O$.
 That is, for each $\vf\in C(\bo)$, there exists a unique $\GG$-harmonic function
 $u\in C(\ob)$ such that 
 \medskip
 
 (i) \ \ $u\bigr|_{\O}$ is $\GG$-harmonic, and
 \medskip
 
  (ii) \ \ $u\bigr|_{\bo} = \vf.$ 
 }

 \medskip
This is the special case Theorems 16.1 of Theorem 13.1  in [\HLse]. 
There are many interesting examples. See [\HLse] for a long list.

Boundary convexity is not required for uniqueness, only an empty core for $X$.
As usual uniqueness is immediate from comparison.

\Theorem{\HH.7. (Comparison)}  {\sl
Suppose that Core$_\GG(X) = \emptyset$ and $K\ss X$ is compact.  If $u\in \PSH_\GG(K)$
and $v\in{\wt \PSH}_\GG(K)$, then the zero maximum principle holds, that is,}
$$
u+v\ \leq\ 0\ \ \ {\sl on}\ \partial K\qquad\Rightarrow \qquad u+v \ \leq \ 0\ \ \ {\sl on}\ K.
\eqno{({\rm ZMP})}
$$
\medskip
\noindent
{\bf Outline of proof.}
By definition $u,v\in \USC(K)$ and on the interior $\Int K$, $u$ is $\GG$-\psh  and 
$v$ is dually $\GG$-plurisubharmonic.  The appropriate notion of {\sl strict} plurisubharmonicity
for general upper semi-continuous functions plays a crucial role, and will be discussed below after outlining
its importance.  If (ZMP) holds for all compact $K\ss X$ under the additional assumption that $u$
is $\GG$-strict, we say that {\sl weak comparison holds for $\GG$ on $X$}.  This weakened version
of comparison has one big advantage, namely that {\sl local implies global} (Theorem 8.3 in [\HLse]).
The proof of completed by showing two things. First,
$$
{\rm Weak\ comparison\ is\ true\ locally}.
\eqno{(\HH.6)}
$$
This follows by a argument based on the ``Theorem on Sums''   -- see Section 10 in [\HLse].
Second, {\sl strict approximation} holds.  That is, since Core$_\GG(X)=\emptyset$, $X$
supports a $C^2$ strictly $\GG$-\psh function $\psi$, and 
$$
\eqalign
{
{\rm If}\ u\ {\rm is\ } &\GG\!-\!{\rm plurisubharmonic}, \cr
{\rm then\ } u+\e\psi\ 
{\rm is\ strictly\ }  &\GG\!-\!{\rm plurisubharmonic, \ for \  each\ } \e>0.
}
\eqno{(\HH.7)}
$$
This follows easily from the definition of  strictness.  Using weak comparison and strict 
approximation, one shows that in the limit comparison holds.\qed

\bigskip
\centerline{\bf Strictness}
\medskip

\Def{\HH.8}  A function $u\in\USC(X)$ is {\bf strictly $\GG$-\psh} if each point in $X$ has a neighborhood
$U$ along with a constant $c>0$ such that for each point $x\in U$ and each test function $\vf$
for $u$ at $x$ 
$$
\tr_W \Hess_x \vf \ \geq\ c \fa W\in \GG_x.
\eqno{(\HH.8)}
$$

To see that this definition of strict agrees with the one given in [\HLse, Def. 7.4,]  one must
compare (\HH.8) with distance in $\Sym(T_xX)$.  For this first note that for $W\in G(p,T_xX)$
the (signed) distance of a point $A\in \Sym(T_xX)$ to the boundary of the positive half-space
defined by the unit normal ${1\over p}P_W$ is simply $\bra A {{1\over p}P_W}$.  Consequently,
the distance from $A\in \cp(\GG_x)$ to $\sim\cp(\GG_x)$ is given by 
$$
\dist(A, \sim\cp(\GG_x)) \ =\ \inf_{W\in \GG_x}  \bra{A}{\smfrac 1 p P_W}
\ =\ \inf_{W\in \GG_x} \smfrac 1 p \tr_W A.
\eqno{(\HH.9)}
$$

For each fixed  $c>0$, $c$-strictness is a subequation.  Therefore, all   the properties in Theorem
\FF.5 hold for $c$-strict $\GG$-\psh functions.  Moreover, as noted in Lemma \HH.3, if $\vf$
is convex and strictly increasing, the composition property holds. Finally, strictness is ``stable''.

\Lemma{\HH.9. ($C^\infty$-Stability Property)} {\sl
Suppose $u$ is strictly $\GG$-\psh and $\psi\in C^\infty(X)$ with
compact support. Then $u+\e\psi$ is  strictly $\GG$-\psh  for all $\e$ sufficiently small.
}
\pf
This is Corollary 7.6 in [\HLse].\qed

 \vskip .3in

\noindent{\headfont \II.\   Geometric Subequations Involving all the Variables.}
\medskip

This is a concept which distinguishes, for example,  the full Laplacian on $\rn$,
which involves all the variables,  from the $p^{\rm th}$ partial Laplacian $\D_p$, which does not.
We shall first treat the euclidean case    (see Section 2 of [\HLf]).
The results will then be carried over to a general riemannian manifold $X$.

Fix a finite dimensional inner product space $V$ and suppose $\GG\ss G(p,V)$ is a closed subset 
of the grassmannian.   
Let Span$\,\GG$ denote the span in $\Sym(V)$ of the elements $P_W$ with $W\in\GG$, and let
$\cp_+(\GG)$ denote the convex cone on $\GG$ with vertex at the origin in $\Sym(V)$.
 Examples show that  $\Span\,\GG$ is often a proper vector subspace of $\Sym(V)$, 
in which case $\cp_+(\GG)$ will have no interior in $\Sym(V)$.  However, considered as a subset of the vector space
$\Span\,\GG$, the interior of $\cp_+(\GG)$ has  closure equal to $\cp_+(\GG)$. We define $\Int_{0} \cp_+(\GG)$ to be
the   interior of $\cp_+(\GG)$ in $\Span\,\GG$  (not in $\Sym(V)$).  In particular, $\Int_{0} \cp_+(\GG)$ is never empty,
and $\cp_+(\GG) = \overline{\Int_{0} \cp_+(\GG)}$.

 By Definition \BB.1, $\cp(\GG) = \{B\in\Sym(V) : \bra B {P_W} \geq0 \ {\rm for\ all\ } W\in\GG\}$.
 Hence, $\cp(\GG)\ss H(A)$ for each closed half-space 
 $H(A) \equiv\{B\in \Sym(V) : \bra AB \geq 0\}$ determined by a non-zero $A\in \cp_+(\GG)$.
 This proves that
 $$
 \cp(\GG)\ =\ \bigcap_{A\in  \cp_+(\GG)} H(A),
 $$
 i.e.,  $\cp(\GG)$ is the ``polar'' of $\cp_+(\GG)$. 
 (Therefore, by the Hahn-Banach/Bipolar Theorem $\cp_+(\GG)$ is the polar of $\cp(\GG)$.)

Since $\cp_+(\GG) = \overline{\Int_0 \cp_+(\GG)}$, this intersection can be taken 
over the smaller set of $A\in \Int_0 \cp_+(\GG)$.
That is,
$$
\cp(\GG)\ =\ \bigcap_{A\in \Int_0 \cp_+(\GG)} H(A).
\eqno{(\II.1)}
$$
This is what will be used below, since
the involvement of all the variables in $\GG$ insures that such $A$ are positive definite, i.e., the linear 
 operators $\bra A{D^2 u}$ are uniformly elliptic.

 The linear operator $\D_A u \equiv \bra A {D^2u}$ with $A\geq 0$ will be referred to as the 
 {\bf $A$-Laplacian}.  Note that from our set theoretic point of view, the subequation
 $\D_A \ss \Sym(V)$ is precisely the closed half-space $H(A)$.

 The following is a restatement of  Proposition 2.8 in  [\HLf] (see also Remark 4.8, page 874 of [K]).

 \Lemma{\II.1} {\sl
 The following are equivalent ways of defining the concept that {\bf $\GG$ involves all the variables}.\medskip
 
 \item{(1)} The only vector $v\in\Sym(V)$ with $v\perp W$ for all $W\in\GG$ is $v=0$.
 \smallskip 
 
 \item{(2)} For each unit vector $e\in V$, $P_e$ is never orthogonal to $\Span \GG$.
 \smallskip 
 
 \item{(3)} There does not exist a hyperplane $W\ss V$ with $\GG\ss\Sym(W)\ss \Sym(V)$.
 \smallskip

 \item{(4)} $\Int_0 \cp_+(\GG) \ss \Int \cp$, i.e., each $A\in \Int_0 \cp_+(\GG)$ is positive definite.
 \smallskip

 \item{(5)} There exists $A\in \Span\,\GG$ with $A>0$.
 }

\medskip

In Section 2 of  [\HLf] such subsets $\GG$ were called ``elliptic''.

We shall apply Lemma \II.1 to the case $V= T_xX$ on a riemannian manifold $X$.
We say that a closed subset $\GG\ss G(p,TX)$
{\bf involves all the variables} if each fibre $\GG_x\ss G(p,T_xX)$ involves all the variables in  the 
vector space $V\equiv T_xX$.
For any smooth section $A(x)\geq0$ of $\Sym(TX)$ the linear operator 
$$
\D_A u \ \equiv\ \bra {A(x)}{\Hess_x u}
$$
will again be referred to as the {\bf $A$-Laplacian}.

Recall  from Definition \FF.8  that $\GG$ {\bf is regular} if each element 
$W_0\in  \GG_{x}$ can be locally 
extended to a smooth section $W(y)$ of $\GG$.  This immediately
implies that  {\sl  each element 
$A_0\in   \cp_+(\GG_{x})$ can be locally 
extended to a smooth section $A(y)$ with $A(y)\in  \cp_+(\GG_y)$},
(since $A_0=\sum_k t_k W_k$ for $t_k>0$ and $W_k \in\GG_{x}$).
Furthermore, if $A(x)>0$, then $A(y)>0$ for $y$ near $x$.
This proves the following.

\Lemma{\II.2} {\sl  
Suppose $\GG\ss G(p,TX)$ is a closed subset involving all the variables and that
$\GG$ is regular.  Then
$$
\cp(\GG_x)\ =\ \bigcap H(A(x)) \qquad{\sl for\ each\ } x\in X
\eqno{(\II.1)'}
$$
where the intersection is taken over all smooth $\cp_+(\GG)$-valued section
$A(y)$ where $A(y)>0$ for $y$ near $x$.}

\Cor{\II.3}  {\sl
A function $u\in\USC(X)$ is $\GG$-\psh  \ $\iff$\  $u$ is $\D_A$-subharmonic 
for each smooth (local) section $A$ of  $\Sym(TX)$  with values in $\cp_+(\GG)$
and $A>0$.
}

\pf If $A$ is a section of $\cp_+(\GG)$, then 
$\cp_+(\GG)\ss \D_A$ over a neighborhood $U$ of $x$, 
so that each $\GG$-\psh function on $U$ is automatically
$\D_A$-subharmonic.  Conversely, if $u$ is $\D_A$-subharmonic for each (local) section $A$ of 
 $\cp_+(\GG)$ with $A>0$, and if $\vf$ is a test function for $u$ at $x$, then
 $\Hess_x \vf \in H(A(x))$, and therefore by (\II.1)$'$, 
$\Hess_x \vf \in \cp(\GG_x)$.\qed

\Note{\II.4}  The simple argument just given also shows the following.  {\sl Suppose $F$ is a subequation
on $X$ which can be written as an intersection of subequations $F=\bigcap_\a F_\a$. 
Then for $u\in\USC(X)$, $u$ is $F$-subharmonic if and only if $u$ is $F_\a$-subharmonic for all $\a$.}
\medskip

Corollary \II.3 has many consequences.  We mention one.

\Theorem{\II.5. (The Strong Maximum Principle for $\GG$-Plurisubharmonic Functions)}  {\sl
Suppose $\GG\ss G(p,TX)$ is regular and involves all the variables.
Then for any compact subset $K$ with $\Int K$ connected and $K=\overline{\Int K}$,
if $u\in \PSH_\GG(K)$ has an interior maximum point, then $u\bigr|_K$ is constant.
}

\pf
Unlike the maximum principle, if the strong maximum principle is true locally,
it is true globally.  However, locally we have $\cp(\GG) \ss \D_A$ with $A>0$,
so the (SMP) for $\D_A$ implies the (SMP) for $\cp(\GG)$.\qed
\medskip

We provide  an example which shows that if the core is non-empty and 
the equation does not involve all the variables, then the (MP), and hence the 
(SMP) can  fail.

\Ex{\II.6}
Let $X\ss \bbr^{n+1}$ be the unit sphere
 $S^n=\{ (x_1,...,x_n, y) \in \rn\times \bbr : x_1^2+\cdots+x_n^2+y^2=1\}$ with the points $y=\pm1$ removed.
Let
 $
H \ =\ \ker \left (dy\bigr|_{TX}\right)
 $
 be the field of ``horizontal''  $(n-1)$-planes on $X$ tangent to the foliation by the 
latitudinal spheres $\{y={\rm constant}\}$, 
and set $\GG_z = \{H(z)\}$ for $z\in S^n$ so that $\GG\ss G(n-1,TX)$.
Calculation shows that for a smooth function $\vf$ defined in a neighborhood
of $X$, 
$$
(\Hess^X \vf)(V,W) \ =\ (\Hess^{\bbr^{n+1}} \vf)(V,W) - \bra VW \nu\cdot \vf
$$
where $\nu$ is the outward-pointing unit normal to $X$.

Now let $\vf =\half (1-y^2)$. Then for $V,W\in H(z)$ horizontal vector fields, the first term vanishes and
the second term yields 
$$
(\Hess^X \vf)(V,W) \ =\    y^2 \bra{V}{W}
$$
Hence $\tr_W \{\Hess^X \vf\} = (n-1)y^2$, proving that $\vf \in\PSH^{\infty}_{\GG}(X)$ and that it
is $\GG$-strict outside $y=0$.
Therefore,  {\sl the maximum principle fails for $\GG$-\psh functions on any domain $\O\ss\ss X$
which contains $S_0^{n-1} \equiv \{y=0\}$ in its interior.}
For any such domain, 
$$
S_0^{n-1} \ss \Core(\O)
$$
because $S^{n-1}_0$ is  a compact minimal $\GG$-submanifold 
and therefore  
any $\GG$-\psh function restricted to it must be constant. (See Theorem \FF.9.)

Note that   $\tr_H \{ \Hess^X u\} \geq0$ is a linear subequation of constant rank
and therefore locally jet equivalent to the partial Laplacian $\D_n$ in local coordinates (Proposition B.3).
Consequently, this subequation satisfies weak local comparison (see the discussion of the proof of Theorem \HH.7).  However, it does 
not satisfy comparison since it does not satisfy the maximum principle.  

 We note that the maximum principle also fails for the subequation consisting of all the 
 $p$-dimensional linear subspaces of $\GG$ (given above), for any $p$, $1\leq p \leq n-1$.

 \vskip .3in

\noindent{\headfont \JJ.\   Distributionally   $\GG$-Plurisubharmonic Functions.}
\medskip

It is easy to see that for the $p^{\rm th}$ partial Laplacian $\D_p$ on $V=\rn$, $p<n$,  there are lots of 
distributional subharmonics (i.e.,  distributions $u$ with $\D_p u$ a non-negative measure) which are not
upper semi-continuous, and hence cannot be horizontally subharmonic.  However, if a closed set $\GG\ss G(p,V)$
involves all the variables, then the appropriate distributional definition of $\GG$-plurisubharmonicity, although
technically not equal, is equivalent to Definition \FF.1. This constant coefficient 
result was proved in Corollary 5.4 of [\HLf].
In this section we extend the result to the variable coefficient case.

First we give the distributional definition.

\Def{\JJ.1}  A distribution $u\in\cd'(X)$ on a riemannian manifold $X$ is
{\bf distributionally $\GG$-\psh} if $\D_A u \geq0$ (a non-negative measure) for every smooth
 section $A(x)$ of $\Sym(TX)$ taking values in $\cp_+(\GG)$.

\medskip 

This distributional notion can not be the ``same'' as $\GG$-plurisubharmonicity, but it is equivalent in a sense we now make
precise.  We exclude the $\GG$-\psh functions which are $\equiv -\infty$ on any component of $X$.  
Let $\lloc(X)$ denote the space of locally integrable functions on $X$.

\Theorem{\JJ.2} {\sl
Assume that $\GG\ss G(p,TX)$ involves all the variables and is regular.
\medskip

(a)  \ \ Suppose $u$ is $\GG$-\psh. Then $u\in \lloc(X)\ss\cd'(X)$, 
and $u$ is distributionally $\GG$-plurisubharmonic.
\medskip

(b)  \ \ Suppose $v\in \cd'(X)$ is distributionally $\GG$-plurisubharmonic.
 Then $v\in \lloc(X)$, and  there exists a unique upper semi-continuous representative
 $u$ of the $\lloc(X)$-class $v$ which is  $\GG$-plurisubharmonic. In fact,
 $$
 u(x) \ = \ {\rm ess}\, \limsup_{y\to x} v(x)
 $$
 is actually independent of $\GG$.
}
\pf Under the hypothesis of Theorem \JJ.2 we can use the next  proposition along with Corollary \II.3
 to reduce to proving 
the  analogous result  for $A$-Laplacians $\D_A$ where $A(x)$ is a smooth  section 
 of  $\Sym(TX)$  having the additional property that $A(x)>0$, i.e., $\D_A$ is uniformly elliptic.

\Prop{\JJ.3}  {\sl
A distribution  $v\in \cd'(X)$ is distributionally $\GG$-\psh  \ $\iff$\  $v$ is distributionally 
$\D_A$-subharmonic for each  smooth (local) section $A$ of  $\Sym(TX)$  with values in $\cp_+(\GG)$
and $A>0$.}

\pf
Suppose $A$ is a local smooth section of $\cp_+(\GG)$ with $A(y)>0$ as in Definition \II.1.
Fix  $x\in X$ and note that  since  $\GG_{x}$ involves all the variables,
 there exists $S_0 \in \Int_0\cp_+(\GG_{x})$ and $S_0>0$ (Lemma \II.1 (4)).
By  the regularity of $\GG$  there exists a local section $S(y)$ of $\cp_+(\GG)$ extending 
$S_0$. Since $S_0>0$, we have that $S(y)>0$ in a neighborhood $U$ of $x$.
Now  for each $\e>0$, $(A+\e S)(y) >0$ on $U$.  That is, locally any smooth section
taking values in $\cp_+(\GG)$ can be approximated by $\cp_+(\GG)$-valued sections 
which are positive definite. Assuming $\D_{A+\e S} u \geq 0$, this implies $\D_A u\geq 0$.\qed

\medskip
\noindent
{\bf   Completion of the Proof of  Theorem  \JJ.2.}
First note that this is a local result.  Note that for  each positive definite 
 $\cp_+(\GG)$-valued  section  $A(x)$,   the $A$-Laplacian
$\D_A$ is of the form
$$
\D_A u \ =\ a(x)\cdot D_x^2u + b(x)\cdot D_xu
$$
where $a(x)$ is a positive definite $n\times n$ matrix and $b(x)$ is $\rn$-valued.  Now the analogue of Theorem
\JJ.2,  with $\GG$-plurisubharmonicity replaced by $\D_A$-subharmonicity, is true.
Details can be found     in  Appendix A of  [\HLs].
An important point in the proof of Theorem \JJ.2 (b) is that the upper semi-continuous 
representative $u$ provided by Appendix A in [\HLs] for a $\D_A$-subharmonic 
distribution $v$ is the same for all sections  $A(x)>0$, since it is the ess-limsup regularization
of the $\lloc$-class $v$. \qed
\medskip
 
 \Remark{\JJ.4} 
The   $\D_A$-harmonics are smooth,  and the notion of $\D_A$-subharmonicity is also equivalent to the
self-defining  notion ``sub-the-$\D_A$-harmonics''  -- again see Appendix A in  [\HLs].
\medskip

 The following gives an easily verified criterion for the regularity of $\GG$.
\medskip\noindent
{\bf Exercise \JJ.5.}  {\sl
 Suppose $\GG\ss G(p,TX)$ is a closed subset which is a smooth fibre-wise neighborhood retract
 in $G(p,TX)$.  Then $\GG$ is regular.
 }
 \medskip

Also note that $\GG$ is a  smooth fibre-wise neighborhood retract in $G(p,TX)$ if and only if it is
a  smooth fibre-wise neighborhood retract in $\Sym(TX)$.

\bigskip
\centerline{\bf Strictness}
\medskip

Recall  that $\GG$-strictness for $u\in\USC(X)$ was defined in Section \HH.  The requirement was that 
locally there exists $c>0$ with $u$ $c$-strict as defined by (\HH.8).
Corollary \II.3 extends to $c$-strictness as follows.

\Prop{\JJ.6} {\sl
A function $u\in\USC(X)$ is $c$-strictly $\GG$-\psh \ \ $\iff$\ \ $u$ is a $c$-strict
$\D_A$-subharmonic  function for each smooth (local) section $A$ of $\cp_+(\GG)$
with $A>0$ at each point.
}

\medskip 

By {\bf $u$ is $c$-strict for $\D_A$}  we  mean that at each point $x$ and for each viscosity
test function $\vf$ for $u$ at $x$, we have $(\D_A \vf)(x) \geq c$.

A distribution $v\in \cd'(X)$ is said to be {\bf  $c$-strict for $\D_A$} (an $A\geq0$ Laplacian) if
$$
\D_A v \ \geq\ c\ \ \ {\rm ( as\ an\ inequality\ of\ measures)}.
\eqno{(\JJ.1)}
$$

If this inequality is true for every smooth section $A$  of $\cp_+(\GG)$, then
$v$ is {\bf $c$-strict as a $\GG$-\psh distribution}.
Proposition \JJ.3 easily extends to 

\Prop{ \JJ.7} {\sl
A distribution $v\in \cd'(X)$ is $c$-strict  for $\GG$ \ \ $\iff$\ \ 
$v$ is $c$-strict  for $\D_A$ for each smooth section $A$ of $\cp_+(\GG)$ which is positive definite.}
\medskip

Since $c$-strictness for the $A$-Laplaican, when $A$ is positive definite, can be show to be
equivalent whether interpreted with viscosity test functions or distributional test functions,
Theorem \JJ.2 has a obvious extension to the $c$-strict  case ($c>0$).
The remainder of the proof is left to the reader, but here is the statement.

\Theorem{\JJ.8}  {\sl
In either part (a) or part (b) of Theorem \JJ.2, if the function in the hypothesis is assumed
to be  $c$-strict, one has $c$-strictness in the conclusion.}

\medskip
Finally we state a result, due to Richberg [R] in the complex case, which carries over to
the $\GG$-\psh case, assuming the following local approximation is possible.

\Def{\JJ.9}  We say that $\GG$ has the {\bf local $C^\infty$-approximation property}
if each point $x\in X$ has a neighborhood $U$ such that for all $u\in C(U)\cap \PSH_\GG(U)$
which are $c$-strict, and all compact $K\ss U$ and $\e>0$, there exists $\wt u\in \PSH^{\infty}_{\GG}(U)$
which is $c$-strict, with 
$
u\leq \wt u \leq u+ \e 
$
on $K$.

\Theorem{\JJ.10} {\sl
Suppose $\GG$ has the local $C^\infty$ strict approximation property, and let $c,\e \in C(X)$
be any given continuous functions satisfying $c>0$ and $1>\e>0$ on $X$.  If $u\in C(X)\cap \PSH_\GG(X)$
is $c$-strict, then there exists $\wt u\in \PSH^{\infty}_{\GG}(X)$, which is $(1-\e)c$-strict, with
$$
u\leq \wt u \leq u+\e \qquad{\rm on\ \ } X.
$$
}

\medskip
The proof in Chapter I, Section 5 of [D], given in the complex case, carries over to this 
much more general case. (See also [GW].)

\vskip .3in

\centerline{\bf Appendix A.   Geometric Subequations}\bigskip

Let $X$ be a riemannian manifold and 
consider a closed subset 
$$
\GG\ \ss\ G(p,TX)
$$
of the Grassmannian of tangent $p$-planes.  The natural candidate for a subequation
$F=F(\GG)$ associated with $\GG$ is defined by its fibres
$$
F_x \ =\ \{ A\in\Sym(T_xX) : \tr_W A\geq 0 \ \ \forall\, W\in \GG_x\}.
\eqno{(A.1)}
$$

For each $W\in \GG_x$ the  condition $\tr_W A\geq0$ defines a closed half-space
(with boundary a hyperplane through 
the origin).  Consequently,
$$
F_x \ {\rm is\  a\  closed\ cone\  with\ vertex\  at \ the\ origin,\ and}
\eqno{(A.2)}
$$
$$
\Int F_x \ =\ \{ A\in\Sym(T_xX) : \tr_W A > 0 \ \ \forall\, W\in \GG_x\}.
\eqno{(A.3)}
$$

Let $\cp_x$ denote the set of non-negative elements in $\Sym(T_xX)$.  Since $\tr_W P\geq0$
for all $W\in G(p, T_xX)$ when $P\in \cp_x$, the fibres $F_x$ defined by (A.1)
satisfy the important {\bf positivity condition}
$$
F_x+\cp_x\ \ss\ F_x.
\eqno{(P)}
$$

Therefore the  fibres $F_x$ satisfy all of the properties of a   constant coefficient (euclidean) 
pure second-order subequation.

\Prop{A.1}{\sl
\medskip

(1)\ \ \  $F_x + \Int \cp_x =  \Int F_x$

\medskip

(2)\ \ \  $F_x =\overline{\Int F_x}$ 

\medskip

(3)\ \ \  $\Int F_x +  \cp_x =  \Int F_x$

\medskip

(4)\ \ \  $A \in \Int F_x  \quad\iff \quad$ there exists a neighborhood of $A$ in $F_x$ of the form

\hskip 1.9in  $\cn_\e(A) \equiv A-\e I +\Int \cp_x$ for some $\e>0$.

}

\pf
(4) \  Note that $\cn_\e(A)$ is an open set containing $A$, and that if 
$A-\e I\in F_x$, then the positivity condition (P) implies that $\cn_\e(A)\ss F_x$.
 \medskip

(1)\ By positivity $F_x+\Int\cp_x \ss F_x$, and it is open since it is the union over $A\in F_x$ of open sets.
Hence it is contained in $\Int F_x$.
Finally, $\Int F_x \ss F_x +\Int \cp_x$ by (4).

\medskip

(2)\ If $A\in F_x$, then by (1) we have $A+\e I \in \Int F_x$ for all  $\e>0$.
Hence, $A = \lim_{\e\to\infty}(A+\e I) \in \overline{\Int F_x}$ proving that $F_x\ss 
\overline{\Int F_x}$.  Since $F_x$ is closed,  we have equality.

\medskip

(3)  \ The containment ``$\ss$'' is proved as in the first half of (1).
The containment ``$\supset$'' follows from $0\in \cp_x$.
\qed

\medskip

Recall the following definition from [\HLse].

\Def{A.2}  A (general)  subset $F\ss \Sym(TX)$ is called a {\bf subequation}
if it satisfies the {\sl positivity condition:}
$$
F_x+\cp_x\ \ss\ F_x \fa x\in X
\eqno{(P)}
$$
and the {\sl three topological conditions:}
$$
{\rm (T_1)}\quad F\ =\ \overline{\Int F}, \qquad 
{\rm (T_2)}\quad F_x\ =\ \overline{\Int F_x}, \qquad
{\rm (T_3)}\quad \Int F_x\ =\  ({\Int F})_x.
$$
(Here $\Int F_x$ means the interior relative to the fibre $\Sym(T_xX)$.)
\medskip

Although $F(\GG)$ is not always  closed (See Proposition  A.6), we shall see
that conditions (T$_1$) and (T$_2$) are always true.  They will be a consequence  of the 
following half of (T$_3$).

\Lemma{A.3} {\sl
The condition
$$
{\rm (T_3)}'\quad \Int F_x\ \ss\  ({\Int F})_x
$$
holds for any closed subset $\GG\ss G(p,TX)$.
Consequently, if a smooth function is $\GG$-strict at a point, 
then it is $\GG$ strict in a neighborhood of that point.
}
\medskip
The proof is given at the end of this appendix.

\Cor{A.4} {\sl The set $F=F(\GG)$ satisfies
$$
{\rm (T_1)}' \quad F\ \ss\ \overline{\Int F}, \qquad 
{\rm (T_2)}\quad F_x\ =\ \overline{\Int F_x}, \qquad
{\rm (T_3)}\quad \Int F_x\ =\  ({\Int F})_x.
$$
}
\pf
Condition  (T$_3)'$ implies  (T$_3$) since $(\Int F)_x$ is an open subset of $F_x$, and hence 
contained in $\Int F_x$.  Property  (T$_2$) is just condition (2) in Proposition A.1.
Finally, by  (T$_2$) and  (T$_3$) we have 
$F_x = \overline{\Int F_x}  = \overline{(\Int F)_x} \ss   \overline{\Int F}$ which
proves  (T$_1)'$.
\qed

 \medskip
We can characterize the case where $F(\GG)$ is closed.

\Def{A.5} The restricted projection $\pi: \GG \to X$ is a {\bf local surjection} if for each $W\in \GG$
and each neighborhood $U$ of $W$, the image $\pi(U\cap \GG)$ contains a neighborhood of $\pi(W)$.
In this case we say that $\GG$ has the {\sl local surjection property}.

\Prop{A.6} {\sl \ \ $F(\GG)$ is closed 
\ \ $\iff$\  \ $\pi:\GG\to X$ is a local surjection.
}

The proof is given at the end of this appendix.

\Cor{A.7} {\sl
A closed subset $\GG\ss G(p,TX)$ determines a subequation $F(\GG)$ via (A.1)
if and only if $\pi: \GG \to X$ is a local surjection.
}

Consequently, we adopt the following definition.

\Def{A.8}  A subset $F\ss \Sym(TX)$ is a {\bf geometrically determined subequation}
if $F= F(\GG)$ with $\GG$ a closed subset of $\Sym(TX)$ having the local surjection property.

\bigskip
\centerline{\bf Strictness}
\medskip

The concept of strictness given in Definition \HH.8 plays an important role for upper semi-continuous functions,
not just smooth functions  (see Definition (\BB.1)$'$) where the notion is unambiguous.

\Def{A.9. ($c$-Strict)}  For each $c>0$ define $F^c = F^c(\GG)$ to be the subset of $\Sym(TX)$  with fibres
$$
F_x^c \ \equiv\ \{ A\in \Sym(T_xX) : \tr_W A \geq c\  \ \forall \, W\in \GG_x\}.
\eqno{(A.4)}
$$

The identity $I$ is a well defined smooth section of $\Sym(TX)$, and $\tr_W I = p$ 
for all $W\in G(p,TX)$. Therefore,
$$
F^c \  = \ F + {c\over p}\cdot I \quad {\rm (fibrewise \ sum)}.
\eqno{(A.5)}
$$
Consequently, all of the previous results for $F$ remain true for $F^c$ ($c\geq 0$). In particular, we have:

\Theorem{A.10} {\sl
If $\GG\ss G(p,TX)$ is a closed subset with the local surjection property,
then for each $c\geq0$ the set $F^c(\GG)$ is a subequation.
}

\bigskip
\centerline{\bf Proofs}\medskip

\noindent
{\bf Proof of Lemma A.3.} 
Assume we are working in a local trivialization $\Sym(T^*V) \cong V\times \Symn$
over an open subset $V\ss X$ containing $x$.
Then each $A\in \Sym(T_xV)$ determines  a smooth section (also denoted $A$) over $V$.
It suffices to prove the following two claims.

\medskip\noindent
{\bf Claim 1:} \ Given $A\in \Sym(T_xV)$, there exists $c>0$ such that
$$
A\in \Int F_x\quad\Rightarrow\quad A\in  F^c_y\ \ {\rm for\ } y \ {\rm near\ } x.
$$
\pf 
If not, there exist  sequences $\{y_j\}$ in $U$ and $W_j\in \GG_{y_j}$ such that 
$$
\lim_{j\to\infty} y_j\ =\ x \and  \lim_{j\to\infty}  \tr_{W_j} A \ =\  0.
$$
By compactness we can assume that $W_j\to W \in \GG_x$, and by continuity this gives
$\tr_W A = 0$, contradicting our assumption that $A\in \Int F_x$ (see (A.3)).
\qed

\medskip\noindent
{\bf Claim 2:} \ If $A$ is a continuous section of $\Sym(TV)$ 
and if for some $c>0$, $A(y)\in  F^c_y$ for all $y$ near $x$, then 
$A(x)\in\Int F$.
\pf
Since $A(y)\in F_y^c$, setting $\e = {c\over p}$, we have that 
$B(y)\equiv A(y) -\e I \in  F_y$ for all $y$ near $x$. 
The set $\cn \equiv B+\Int \cp$,
defined using fibre-wise sum, is the translation of the open subset $\Int \cp$ of $V\times \Symn$ by a continuous
section.
 Hence, $\cn$ is open in $\Sym(TV)$.  Since $B(y)\in F_y$ for all $y$, we have $\cn \equiv B+\Int \cp \ss F$.
 Hence, $\cn \ss\Int F$. Finally, $A(x) = B(x) +\e I \in \cn$  by positivity.\qed

 \medskip
\noindent
{\bf Proof of  Proposition A.6.}  The  assertion is local so
we may assume that $X$ is an open subset of $\rn$ and $\pi:X\times G(p,\rn)\to X$ is projection onto the
first factor, with $G(p,\rn)\ss \Sym(\rn)$.

Suppose $\pi\bigr|_\GG$ is locally surjective.  Let  $(x_j, A_j)\in F$ be a convergent sequence,
 $x_j\to x, A_j\to A$.
Fix $W\in \GG_x$.  By hypothesis for each neighborhood $N_\d(W)$ of $W$,   $\pi\{(X\times N_\d(W))\cap\GG\}$ contains a neighborhood of $x$.
Hence we may pick $W_j\in \GG_{x_j}$ with $W_j\to W$.  Since $\tr_{W_j}A_j \geq 0$ for all $j$ we
have  $\tr_W A\geq0$, and so $A\in F_x$.

For the converse, suppose $\pi\bigr|_\GG$ is not locally surjective.  Then there exists $(x,W)\in \GG$ 
and a neighborhood $N(W)$ of $W$   in  $G(p,\rn)$  so that 
$\pi\{(X\times N(W))\cap\GG\}$ does not contain a neighborhood of $x$.
Hence there exists  a sequence of points  $x_j\to x$ in $X$, 
such that $\GG_{x_j} \cap N(W) =\emptyset$ for all $j$.

If $\e>0$ is chosen small enough, then for all $V\in G(p,\rn)$
$$
\bra {P_V}{P_{W^\perp}} \ <\   \e p  \qquad\Rightarrow\qquad V\in N(W).
\eqno{(A.1)}
$$
Consequently, since $\bra {P_V}{ -P_W} \geq -p$, we have that
$$
 V\notin N(W)    \qquad\Rightarrow\qquad  \bra {P_V}{-P_W +\smfrac 1 \e  P_{W^\perp}} \ \geq \ 0.
\eqno{(A.2)}
$$
Since $\GG_{x_j} \cap N(W) = \emptyset$, this proves that $A\equiv -P_W+{1\over \e} P_{W^\perp} \in F_{x_j}$.
However, $\bra A {P_W} = -1$ and $W\in \GG_x$, and so $A\notin F_x$.  We conclude that  $F$ is not closed.\qed

 \vfill\eject

\centerline{\bf Appendix B.   The Linear Geometric Case.}\bigskip

In this appendix we consider the extreme geometric case where each $\GG_x$ is a single
point $W_x\in G(p,T_xX)$, or equivalently, each $\cp(\GG_x)$ is the half space in $\Sym(T_xX)$
with inward normal $P_{W_x}$ (orthogonal projection onto $W_x$).  Said differently, the subequation
$\cp(\GG)$ is linear and given by the {\bf $W$-Laplacian}
$$
\left(  \D_W u  \right)(x)\ =\  { \bra {P_{W_x}} {\Hess_x u}}_{\rm riem.}
\ =\ \tr_{W_x} \Hess_x u.
\eqno{(B.1)}
$$
It is more appropriate to refer to {\bf $W$-subharmonic} functions, rather than
$\GG$-\psh functions in this {\bf linear-geometric case}.

\Ex{B.1.  (The $p^{\rm th}$ Horizontal Laplacian)} 
In this example, choose a single $p$-plane $W\in G(p,\rn)$, which might as well be the first
coordinate $p$-plane $W\equiv \bbr^p\times \{0\}\ss \rn$.  Abbreviate $P_{\bbr^p\times \{0\}}$ to $P$.
Then
$$
\D_P u\ =\ \bra P {D^2u}\ =\ \tr_P D^2 u\ =\ \sum_{j=1}^p {\partial^2 u\over \partial x_j^2}
\eqno{(B.2)}
$$
is the {\bf $p^{\rm th}$ horizontal Laplacian}.  The terminology ``horizontally subharmonic'' and ``horizontally $p$-convex'' is appropriate in this case.
\medskip

Suppose $h$ and $H$ are smooth functions defined on an open subset on $\rn$, 
with $h$ taking values in GL$_n(\bbr)$ and with $H$ taking values in Hom$(\rn, \Sym(\rn))$.

\Def{B.2}  An equation of the form
$$
L u\ =\ \bra {h^t P h} {D^2u} + \bra{H^t(P)}{Du}
\eqno{(B.3)}
$$
is said to be {\bf jet equivalent to} $\D_p$.\medskip

The linear-geometric case is jet equivalent to $\D_p$ in any local coordinate system.

\Prop{B.3}  {\sl
If $W$ is  a smooth section of the Grassmann bundle $G(p,TX)$ over $X$, then
the $W$-Laplacian is jet equivalent to the $p^{\rm th}$ horizontal Laplacian over any
local coordinate chart.
}
\pf
Choose a local orthonormal frame field $e_1,...,e_n$ for $\rn$ with $e_1,...,e_p$ a frame for
$W$.  Define $h(x)$ with values in  GL$_n(\bbr)$ by $e= h {\partial\over \partial x}$.
Then, in the given local coordinates,
$$
\D_W  u\ =\ \bra {h^t P h} {D^2u} -  \bra{\G^t (h^tPh)}{Du}
\eqno{(B.4)}
$$
follows from Proposition 5.5 in [\HLse]. \qed

\medskip
\Prop{B.4}  {\sl
A subequation $L$ is locally jet equivalent to the $p^{\rm th}$ horizontal Laplacian  $\D_p$ if and only if
in any local coordinate system
$$
L  u\ =\ \bra {E} {D^2u} -  \bra{b}{Du}
\eqno{(B.5)}
$$
where $b$ and $E$ are smooth and $E(x)$ has  rank $p$ at each point $x$.
}
\eject

\pf
First note that $E=h^t P h$ in (B.3) has constant rank $p$. Conversely, assume $E$ in
(B.5) has rank $p$ at each point.   Then $E(x)$ has a unique smooth square root $A(x)$ in $\Sym(\rn)$.  Let $B$ denote orthogonal projection onto the  null space of $E$.  Then the 
inverse of $A+B$ conjugates $E$ to $P_W$ where $W^\perp \equiv \ker E$.
Finally it is easy to (locally) conjugate $P_W$ to $P$ and find a smooth section $H$ 
of $\Hom(\rn, \Symn)$ with $H^t (P) = b$.\qed

\vskip .2in



\centerline{\bf References}

\vskip .2in

\noindent
\item{[\Cr]}   M. G. Crandall,  {\sl  Viscosity solutions: a primer},  
pp. 1-43 in ``Viscosity Solutions and Applications''  Ed.'s Dolcetta and Lions, 
SLNM {\bf 1660}, Springer Press, New York, 1997.

 \smallskip

\noindent
\item{[\CIL]}   M. G. Crandall, H. Ishii and P. L. Lions {\sl
User's guide to viscosity solutions of second order partial differential equations},  
Bull. Amer. Math. Soc. (N. S.) {\bf 27} (1992), 1-67.

 \smallskip

\noindent
\item{[D]}  J.-P. Demailly, { Complex Analytic and Differential Geometry},
  e-book at Institut Fourier, UMR  5582 du CNRS,
   Universit\'e de Grenoble I, Saint-Martin d'H\`eres, France:
   can be found at http://www-fourier.ujfgrenoble.fr/~demailly/books.html.
 \smallskip

\noindent
\item{[GW]}  R. E. Greene and H. Wu, { $C^\infty$ approximations of convex, 
subharmonic, and plurisubharmonic functions},
Ann. Scien. De l'E.N.S. 4$^e$ s'erie, tome 12, no. 1 (1979), 47-84.
 \smallskip

\item {[\HLo]} F. R. Harvey and H. B. Lawson, Jr., 
 {\sl  An introduction to potential theory in calibrated geometry}, Amer. J. Math.  {\bf 131} no. 4 (2009), 893-944.   
 
 ArXiv:math.0710.3920.

\smallskip

\item {[\HLt]} \ \----------, {\sl  Duality of positive currents and plurisubharmonic functions in calibrated geometry},  Amer. J. Math.    {\bf 131} no. 5 (2009), 1211-1240.

 ArXiv:math.0710.3921.

\smallskip

\item {[\HLth]}   \ \----------, {\sl  Dirichlet duality and the non-linear Dirichlet problem},    Comm. on Pure and Applied Math. {\bf 62} (2009), 396-443.

\smallskip

\item {[\HLf]}  \ \----------,  {\sl  Plurisubharmonicity in a general geometric context},  Geometry and Analysis {\bf 1} (2010), 363-401. ArXiv:0804.1316.

\smallskip

\item {[\HLfi]}  \ \----------, {\sl  Lagrangian plurisubharmonicity and convexity},  Stony Brook Preprint (2009).
   
 \smallskip

\item {[\HLs]}  \ \----------,  {\sl  Potential Theory on almost complex manifolds}, 
ArXiv: 1107.2584.
\smallskip

\item {[\HLse]}  \ \----------, {\sl Dirichlet Duality and the Nonlinear Dirichlet Problem on Riemannian Manifolds},  
J. Diff. Geom. (2011).

\smallskip

\item {[\HLe]}  \ \----------, 
{\sl  The restriction theorem for fully nonlinear subequations},   

ArXiv:1101.4850.

\smallskip

\item {[\HLn]}  \ \----------, {\sl  Foundations of  $p$-convexity 
and $p$-plurisubharmonicity in riemannian geometry},  Stony Brook Preprint.

\smallskip

\item {[K]} N. V. Krylov, {\sl  On the general notion of fully nonlinear second-order elliptic equations},
 Trans. A.M.S.  {\bf 347}    (1995), 857-895.

\smallskip

\item {[R]} R. Richberg, {\sl  Stetige streng pseudokonvexe Funktionen},
  Math. Ann. {\bf 175}    (1968), 257-286.

\smallskip

\item {[Ro]} C.  Robles, {\sl  Parallel calibrations and minimal submanifolds},
 to appear in Illinois J. Math. (2011).  ArXiv:0808.2158.

\smallskip

\end